\documentclass[titlepage,11pt]{article}
\usepackage{latexsym}
\usepackage{amsfonts}
\usepackage{amsmath}
\usepackage{tikz}
\usepackage{float}
\usepackage{comment}
\usetikzlibrary{decorations.pathreplacing,calligraphy}
\usetikzlibrary{positioning,arrows.meta}

\newcommand\blackslug{\hbox{\hskip 1pt \vrule width 4pt height 8pt depth 1.5pt
        \hskip 1pt}}
\newcommand\bbox{\hfill \quad \blackslug \bigbreak}

\def\DD{\hbox{-}}
\def\CC{\hbox{-}\cdots\hbox{-}}
\def\LL{,\ldots,}
\def\cupcup{\cup\cdots\cup}

\title{The minimal nonplanar strong digraphs}
\author{Stephen Bartell\\
Princeton University, Princeton, NJ 08544
\\
\\
Paul Seymour\thanks{Supported by AFOSR grant
FA9550-22-1-0234, and NSF grant  DMS-2154169.}\\
Princeton University, Princeton, NJ 08544}

\date{February 13, 2025; revised \today}

\newtheorem{thm}{}[section]

\newcommand{\Proof}{\noindent{\bf Proof.}\ \ }

\begin{document}
\maketitle
\begin{abstract}
Kuratowski's theorem says that the minimal (under subgraph containment) graphs that are not planar are the subdivisions of $K_5$ and of $K_{3,3}$.
Here we study the minimal (under subdigraph containment) strongly-connected 
digraphs 
that are not planar. We also find the minimal strongly-connected non-outerplanar digraphs and the minimal 
strongly-connected non-series-parallel digraphs. 
\end{abstract}

\section{Introduction}

What are the minimal digraphs (under subdigraph containment) that are strong and nonplanar? (``Strong'' means 
strongly-connected.) We will 
prove that for every such digraph, its underlying graph is a subdivision of a graph of one of the seven types shown in Figures \ref{fig:nonplanar} and 
\ref{fig:clam}. (We will define these types more precisely later.)

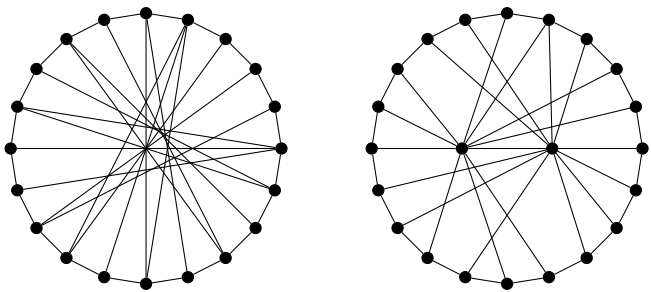
\begin{figure}[H]
\centering

\begin{tikzpicture}[scale=1.2,auto=left]
\tikzstyle{every node}=[inner sep=1.5pt, fill=black,circle,draw]
\begin{scope}[shift ={(-4,0)}]
\def\r{1.5}
\def\a{18}

\node (a1) at ({\r*cos(0*\a)}, {\r*sin(0*\a)}){};
\node (a2) at ({\r*cos(1*\a)}, {\r*sin(1*\a)}){};
\node (a3) at ({\r*cos(2*\a)}, {\r*sin(2*\a)}){};
\node (a4) at ({\r*cos(3*\a)}, {\r*sin(3*\a)}){};
\node (a5) at ({\r*cos(4*\a)}, {\r*sin(4*\a)}){};
\node (a6) at ({\r*cos(5*\a)}, {\r*sin(5*\a)}){};
\node (a7) at ({\r*cos(6*\a)}, {\r*sin(6*\a)}){};
\node (a8) at ({\r*cos(7*\a)}, {\r*sin(7*\a)}){};
\node (a9) at ({\r*cos(8*\a)}, {\r*sin(8*\a)}){};
\node (a10) at ({\r*cos(9*\a)}, {\r*sin(9*\a)}){};
\node (a11) at ({\r*cos(10*\a)}, {\r*sin(10*\a)}){};
\node (a12) at ({\r*cos(11*\a)}, {\r*sin(11*\a)}){};
\node (a13) at ({\r*cos(12*\a)}, {\r*sin(12*\a)}){};
\node (a14) at ({\r*cos(13*\a)}, {\r*sin(13*\a)}){};
\node (a15) at ({\r*cos(14*\a)}, {\r*sin(14*\a)}){};
\node (a16) at ({\r*cos(15*\a)}, {\r*sin(15*\a)}){};
\node (a17) at ({\r*cos(16*\a)}, {\r*sin(16*\a)}){};
\node (a18) at ({\r*cos(17*\a)}, {\r*sin(17*\a)}){};
\node (a19) at ({\r*cos(18*\a)}, {\r*sin(18*\a)}){};
\node (a20) at ({\r*cos(19*\a)}, {\r*sin(19*\a)}){};

\foreach \from/\to in {a1/a2,a2/a3,a3/a4,a4/a5,a5/a6,a6/a7,a7/a8,a8/a9,a9/a10,
a10/a11,a11/a12,a12/a13,a13/a14,a14/a15,a15/a16,a16/a17,a17/a18,a18/a19,a19/a20, a20/a1,a1/a10, a1/a11,a1/a12,a2/a13,a3/a13,a4/a14,
a5/a14,a5/a15, a5/a16,a6/a16,a6/a17,a7/a18,a8/a18,a8/a19, a9/a20,a10/a20}
\draw [-] (\from) -- (\to);
\end{scope}
\def\r{1.5}
\def\a{18}
\node (a1) at ({\r*cos(0*\a)}, {\r*sin(0*\a)}){};
\node (a2) at ({\r*cos(1*\a)}, {\r*sin(1*\a)}){};
\node (a3) at ({\r*cos(2*\a)}, {\r*sin(2*\a)}){};
\node (a4) at ({\r*cos(3*\a)}, {\r*sin(3*\a)}){};
\node (a5) at ({\r*cos(4*\a)}, {\r*sin(4*\a)}){};
\node (a6) at ({\r*cos(5*\a)}, {\r*sin(5*\a)}){};
\node (a7) at ({\r*cos(6*\a)}, {\r*sin(6*\a)}){};
\node (a8) at ({\r*cos(7*\a)}, {\r*sin(7*\a)}){};
\node (a9) at ({\r*cos(8*\a)}, {\r*sin(8*\a)}){};
\node (a10) at ({\r*cos(9*\a)}, {\r*sin(9*\a)}){};
\node (a11) at ({\r*cos(10*\a)}, {\r*sin(10*\a)}){};
\node (a12) at ({\r*cos(11*\a)}, {\r*sin(11*\a)}){};
\node (a13) at ({\r*cos(12*\a)}, {\r*sin(12*\a)}){};
\node (a14) at ({\r*cos(13*\a)}, {\r*sin(13*\a)}){};
\node (a15) at ({\r*cos(14*\a)}, {\r*sin(14*\a)}){};
\node (a16) at ({\r*cos(15*\a)}, {\r*sin(15*\a)}){};
\node (a17) at ({\r*cos(16*\a)}, {\r*sin(16*\a)}){};
\node (a18) at ({\r*cos(17*\a)}, {\r*sin(17*\a)}){};
\node (a19) at ({\r*cos(18*\a)}, {\r*sin(18*\a)}){};
\node (a20) at ({\r*cos(19*\a)}, {\r*sin(19*\a)}){};
\node (c) at (-.5,0){};
\node (b) at (.5,0){};
\foreach \from/\to in {a1/a2,a2/a3,a3/a4,a4/a5,a5/a6,a6/a7,a7/a8,a8/a9,a9/a10,
a10/a11,a11/a12,a12/a13,a13/a14,a14/a15,a15/a16,a16/a17,a17/a18,a18/a19,a19/a20, a20/a1,b/a1,c/a2, c/a3,b/a4,b/a5,c/a5,c/a6,b/a7,
b/a8,c/a9,c/a10,c/a11, b/a12,b/a13, c/a14, b/a15, c/a16, c/a17,b/a18, b/a19,b/a20,c/b}
\draw [-] (\from) -- (\to);

\end{tikzpicture}
\caption{A M\"obius chain and a double wheel.} \label{fig:nonplanar}
\end{figure}
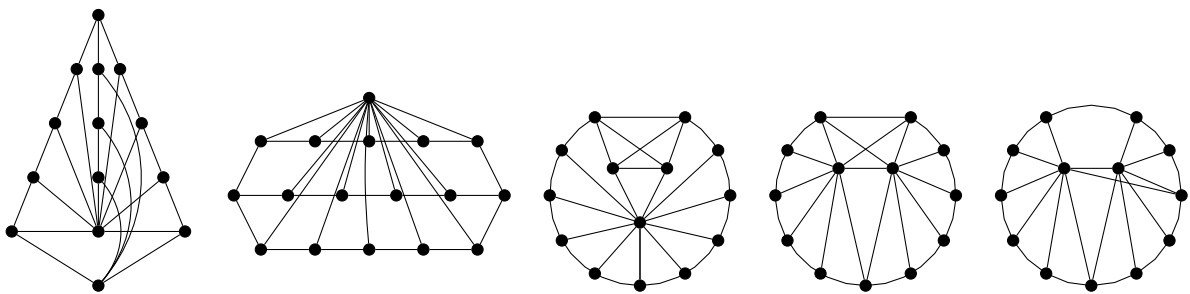
\begin{figure}[H]
\centering

\begin{tikzpicture}[scale=1.2,auto=left]
\tikzstyle{every node}=[inner sep=1.5pt, fill=black,circle,draw]
\begin{scope}[shift={(-3,-2.1)},scale=1.2]

\node (a1) at (0,.5){};
\node (a2) at (0,1){};
\node (a3) at (0,1.5){};
\node (a4) at (0,2){};
\node (a5) at (0,2.5){};
\node (a6) at (0,3){};
\node (a9) at (-.8,1){};
\node (a10) at (-.6,1.5){};
\node (a11) at (-.4,2){};
\node (a12) at (-.2,2.5){};
\node (a13) at (.8, 1){};
\node (a14) at (.6,1.5){};
\node (a15) at (.4, 2){};
\node (a16) at (.2,2.5){};

\foreach \from/\to in {a2/a3,a3/a4,a4/a5,a5/a6,a9/a10,
a10/a11,a11/a12,a13/a14,a14/a15,a15/a16, a12/a6,a16/a6,a1/a9,a1/a13, a2/a9, a2/a10, a2/a11,a2/a12,a2/a13,a2/a14,a2/a15,a2/a16}
\draw [-] (\from) -- (\to);

\draw (a1) to [bend right = 40] (a3);
\draw (a1) to [bend right = 40] (a4);
\draw (a1) to [bend right = 40] (a5);

\end{scope}

\begin{scope}[shift ={(0,-.5)}, scale = .6]
\node (p1) at (-2.5,0){};
\node (p2) at (2.5,0){};
\node (p3) at (0,1.8){};
\node (p4) at (-2,1){};
\node (p5) at (-1,1){};
\node (p6) at (0,1){};
\node (p7) at (1,1){};
\node (p8) at (2,1){};
\node (p9) at (-1.5,0){};
\node (p10) at (-.5,0){};
\node (p11) at (.5,0){};
\node (p12) at (1.5,0){};
\node (p13) at (-2,-1){};
\node (p14) at (-1,-1){};
\node (p15) at (0,-1){};
\node (p16) at (1,-1){};
\node (p17) at (2,-1){};

\foreach \from/\to in {p1/p4,p1/p9,p1/p13,p2/p8,p2/p12,p2/p17,p3/p4,p3/p5,p3/p6,p3/p7,p3/p8,p3/p9,p3/p10,p3/p11,p3/p12,p3/p13,
p3/p14,p3/p16,p3/p17,p4/p5,p5/p6,p6/p7,p7/p8,p9/p10,p10/p11,p11/p12,p13/p14,p14/p15,p15/p16,p16/p17}
\draw [-] (\from) -- (\to);
\draw (p3) to [bend right = 5] (p15);
\end{scope}

\begin{scope}[shift ={(3,-.5)}, scale = 1]
\def\r{1}
\def\a{30}
\node (a3) at ({\r*cos(0*\a)}, {\r*sin(0*\a)}){};
\node (a2) at ({\r*cos(1*\a)}, {\r*sin(1*\a)}){};
\node (a1) at ({\r*cos(2*\a)}, {\r*sin(2*\a)}){};
\node (a11) at ({\r*cos(4*\a)}, {\r*sin(4*\a)}){};
\node (a10) at ({\r*cos(5*\a)}, {\r*sin(5*\a)}){};
\node (a9) at ({\r*cos(6*\a)}, {\r*sin(6*\a)}){};
\node (a8) at ({\r*cos(7*\a)}, {\r*sin(7*\a)}){};
\node (a7) at ({\r*cos(8*\a)}, {\r*sin(8*\a)}){};
\node (a6) at ({\r*cos(9*\a)}, {\r*sin(9*\a)}){};
\node (a5) at ({\r*cos(10*\a)}, {\r*sin(10*\a)}){};
\node (a4) at ({\r*cos(11*\a)}, {\r*sin(11*\a)}){};

\node (a) at (-.3,.3){};
\node (b) at (.3,.3){};
\node (c) at (0,-.3){};

\draw [domain=120:420] plot ({cos(\x)}, {sin(\x)});

\foreach \from/\to in {a/b,a/c,b/c,a1/a11,b/a1,a/a11,c/a2,c/a3,c/a4,c/a5,c/a6,c/a6,c/a7,c/a8,c/a9,c/a10,a/a1,b/a11}
\draw [-] (\from) -- (\to);

\end{scope}

\begin{scope}[shift ={(5.5,-.5)}, scale = 1]
\def\r{1}
\def\a{30}
\node (a3) at ({\r*cos(0*\a)}, {\r*sin(0*\a)}){};
\node (a2) at ({\r*cos(1*\a)}, {\r*sin(1*\a)}){};
\node (a1) at ({\r*cos(2*\a)}, {\r*sin(2*\a)}){};
\node (a11) at ({\r*cos(4*\a)}, {\r*sin(4*\a)}){};
\node (a10) at ({\r*cos(5*\a)}, {\r*sin(5*\a)}){};
\node (a9) at ({\r*cos(6*\a)}, {\r*sin(6*\a)}){};
\node (a8) at ({\r*cos(7*\a)}, {\r*sin(7*\a)}){};
\node (a7) at ({\r*cos(8*\a)}, {\r*sin(8*\a)}){};
\node (a6) at ({\r*cos(9*\a)}, {\r*sin(9*\a)}){};
\node (a5) at ({\r*cos(10*\a)}, {\r*sin(10*\a)}){};
\node (a4) at ({\r*cos(11*\a)}, {\r*sin(11*\a)}){};

\node (a) at (-.3,.3){};
\node (b) at (.3,.3){};

\draw [domain=120:420] plot ({cos(\x)}, {sin(\x)});

\foreach \from/\to in {a/b,a1/a11,b/a1,b/a2,b/a3,b/a4,b/a5,b/a6,a/a6,a/a7,a/a8,a/a9,a/a10,a/a11,a/a1,b/a11}
\draw [-] (\from) -- (\to);
\end{scope}
\begin{scope}[shift = {(8,-.5)}, scale = 1]
\def\r{1}
\def\a{30}
\node (a3) at ({\r*cos(0*\a)}, {\r*sin(0*\a)}){};
\node (a2) at ({\r*cos(1*\a)}, {\r*sin(1*\a)}){};
\node (a1) at ({\r*cos(2*\a)}, {\r*sin(2*\a)}){};
\node (a11) at ({\r*cos(4*\a)}, {\r*sin(4*\a)}){};
\node (a10) at ({\r*cos(5*\a)}, {\r*sin(5*\a)}){};
\node (a9) at ({\r*cos(6*\a)}, {\r*sin(6*\a)}){};
\node (a8) at ({\r*cos(7*\a)}, {\r*sin(7*\a)}){};
\node (a7) at ({\r*cos(8*\a)}, {\r*sin(8*\a)}){};
\node (a6) at ({\r*cos(9*\a)}, {\r*sin(9*\a)}){};
\node (a5) at ({\r*cos(10*\a)}, {\r*sin(10*\a)}){};
\node (a4) at ({\r*cos(11*\a)}, {\r*sin(11*\a)}){};

\node (a) at (-.3,.3){};
\node (b) at (.3,.3){};

\draw [domain=0:360] plot ({cos(\x)}, {sin(\x)});

\foreach \from/\to in {a/b,b/a1,b/a2,b/a3,b/a4,b/a5,b/a6,a/a6,a/a7,a/a8,a/a9,a/a10,a/a11,a/a3}
\draw [-] (\from) -- (\to);

\end{scope}
\end{tikzpicture}
\caption{A conch, a mussel, a scallop, a clam and a whelk.} \label{fig:clam}
\end{figure}

Remarkably, the question reduces, almost exactly, to a question about undirected graphs. Say a nonplanar 
graph $G$ is {\em
almost-planar} if $G$ is 3-connected and $F$ is the edge-set of a forest, where $F$ is the set of edges $e$ such that $G\setminus e$ is nonplanar. 
If $G$ is a digraph that is minimal strong
and nonplanar, and not the directed subdivision of a smaller digraph, then it is an orientation of an almost-planar graph; 
and we will find explicitly
all the almost-planar graphs. (These are the graphs of the two figures.) 
Thus we have:
\begin{thm}\label{planar}
A strong digraph is planar if and only if no subdigraph is a directed subdivision of a strong orientation of either a
M\"obius chain, a double wheel,  a conch, a scallop, a mussel, a clam, or a whelk.
\end{thm}
Let us say a {\em fan} in $G$ is a pair $(P,v)$ where $P$ is a path of $G$ with length at least two, and $v\in V(G)\setminus V(P)$
 is adjacent to every vertex of the interior of $P$ (and possibly also adjacent to its ends), and each vertex in the 
interior of $v$ has degree three in $G$. If $\{v,x,y\}$ is a triangle of a graph, let us subdivide the edge $xy$, replacing it by a path $P$, and 
make $v$ adjacent to all the new vertices introduced (and possibly delete the edges $vx,vy$). Then $(P,v)$ is a fan of the new 
graph, and we call this process {\em opening a fan}. All the five types of graphs shown in Figure \ref{fig:clam} can be obtained 
from graphs of bounded size (in most cases, from $K_5$) 
by opening a fan at most three times, and so in some sense they are ``essentially'' of bounded size; 
but this is not so for the types of Figure \ref{fig:nonplanar}.

\ref{planar} does not quite answer the question of what the minimal strong nonplanar digraphs are. The seven different types of graphs 
in the figures all give rise to infinitely many instances of minimal strong nonplanar digraphs, but not every graph of these
types can be appropriately oriented, and some can be oriented in many ways.  
This is a non-trivial issue, that we discuss in the final section.

This work grew out of the senior thesis~\cite{thesis} of one of the authors, and a much easier question: what are the minimal digraphs that are strong and not outerplanar? 
({\em Outerplanar} means it can be drawn in a plane with all vertices
incident with the infinite region,) For that problem we can list all such digraphs explicitly. Similarly we can 
find explicitly
all minimal digraphs that are strong and not series-parallel (we say a digraph is {\em series-parallel} if it is an orientation of a graph that contains no subdivision of $K_4$ as a subgraph). 

All graphs and digraphs in this paper are assumed to be finite, and to have no loops or parallel edges, and digraphs have no directed cycles of length two.
If $G$ is a digraph and $uv\in E(G)$, {\em subdividing} $uv$ means making a new digraph by deleting the edge $uv$ and adding 
a 
directed path from $u$ to $v$ of new vertices (except for $u,v$), and we call a digraph obtained from $G$ by iterating this process a  {\em directed subdivision} of $G$.

\begin{figure}[H]
\centering

\begin{tikzpicture}[scale=1.5,auto=left]
\tikzstyle{every node}=[inner sep=1.5pt, fill=black,circle,draw]
\node (a1) at (0,2) {};
\node (a2) at (0,0) {};
\node (a3) at (-1.5,1) {};
\node (a4) at (-1, 1) {};
\node (a5) at (-.5, 1) {};
\def\s{1.5}
\draw [-{Stealth[scale=\s]}] (a1) -- (-.75,1.5);
\draw [-{Stealth[scale=\s]}] (a1) -- (-.5,1.5);
\draw [-{Stealth[scale=\s]}] (a5) -- (-.25,1.5);
\draw [-{Stealth[scale=\s]}] (a3) -- (-.75,.5);
\draw [-{Stealth[scale=\s]}] (a4) -- (-.5,.5);
\draw [-{Stealth[scale=\s]}] (a2) -- (-.25,.5);

\foreach \from/\to in {a1/a3,a1/a4,a1/a5,a2/a3,a2/a4,a2/a5}
\draw [-] (\from) -- (\to);

\def\d{2}
\node (a1) at (0+\d,2) {};
\node (a2) at (0+\d,0) {};
\node (a3) at (-1.5+\d,1) {};
\node (a4) at (-1+\d, 1) {};
\node (a5) at (-.5+\d, 1) {};
\def\s{1.5}
\draw [-{Stealth[scale=\s]}] (a1) -- (-.75+\d,1.5);
\draw [-{Stealth[scale=\s]}] (a1) -- (-.5+\d,1.5);
\draw [-{Stealth[scale=\s]}] (a1) -- (-.25+\d,1.5);
\draw [-{Stealth[scale=\s]}] (a3) -- (-.75+\d,.5);
\draw [-{Stealth[scale=\s]}] (a4) -- (-.5+\d,.5);
\draw [-{Stealth[scale=\s]}] (a5) -- (-.25+\d,.5);
\draw [-{Stealth[scale=\s]}] (a2) -- (0+\d,1);

\foreach \from/\to in {a1/a3,a1/a4,a1/a5,a2/a3,a2/a4,a2/a5, a1/a2}
\draw [-] (\from) -- (\to);

\def\e{3}
\node (b1) at  (\e,1) {};
\node (b2) at  (\e+1,0) {};
\node (b3) at  (\e+2,1) {};
\node (b4) at  (\e+1,2) {};

\draw [-{Stealth[scale=\s]}] (b1) -- (\e+.5,.5);
\draw [-{Stealth[scale=\s]}] (b2) -- (\e+1.5, .5);
\draw [-{Stealth[scale=\s]}] (b3) -- (\e+1.5,1.5);
\draw [-{Stealth[scale=\s]}] (b4) -- (\e+.5,1.5);
\draw [-{Stealth[scale=\s]}] (b1) -- (\e+1.3,1);
\draw [-{Stealth[scale=\s]}] (b2) -- (\e+1,1.3);

\foreach \from/\to in {b1/b2,b1/b3,b1/b4,b2/b3,b2/b4,b3/b4}
\draw [-] (\from) -- (\to);
\def\f{7}
\def\r{1}
\def\g{1}
\node(c0) at (\f, \g){};
\node (c1) at ({\f+\r*cos(0)}, {\g+\r*sin(0)}) {};
\node(c2) at ({\f+\r*cos(45)}, {\g+\r*sin(45)}){};
\node(c3) at ({\f+\r*cos(90)}, {\g+\r*sin(90)}){};
\node(c4) at ({\f+\r*cos(135)}, {\g+\r*sin(135)}){};
\node(c5) at ({\f+\r*cos(180)}, {\g+\r*sin(180)}){};
\node(c6) at ({\f+\r*cos(225)}, {\g+\r*sin(225)}){};
\node(c7) at ({\f+\r*cos(270)}, {\g+\r*sin(270)}){};
\node(c8) at ({\f+\r*cos(315)}, {\g+\r*sin(315)}){};

\foreach \from/\to in {c1/c2,c2/c3,c3/c4,c4/c5,c5/c6,c6/c7,c7/c8,c8/c1,c0/c1,c0/c2,c0/c3,c0/c4,c0/c5,c0/c6,c0/c7,c0/c8}
\draw [-] (\from) -- (\to);

\draw [-{Stealth[scale=\s]}] (c1) -- ({\f+\r*(cos(0)+cos(45))/2},{\g+\r*(sin(0)+sin(45))/2});
\draw [-{Stealth[scale=\s]}] (c3) -- ({\f+\r*(cos(90)+cos(135))/2},{\g+\r*(sin(90)+sin(135))/2});
\draw [-{Stealth[scale=\s]}] (c5) -- ({\f+\r*(cos(180)+cos(225))/2},{\g+\r*(sin(180)+sin(225))/2});
\draw [-{Stealth[scale=\s]}] (c7) -- ({\f+\r*(cos(270)+cos(315))/2},{\g+\r*(sin(270)+sin(315))/2});
\draw [-{Stealth[scale=\s]}] (c3) -- ({\f+\r*(cos(90)+cos(45))/2},{\g+\r*(sin(90)+sin(45))/2});
\draw [-{Stealth[scale=\s]}] (c5) -- ({\f+\r*(cos(180)+cos(135))/2},{\g+\r*(sin(180)+sin(135))/2});
\draw [-{Stealth[scale=\s]}] (c7) -- ({\f+\r*(cos(270)+cos(225))/2},{\g+\r*(sin(270)+sin(225))/2});
\draw [-{Stealth[scale=\s]}] (c1) -- ({\f+\r*(cos(0)+cos(315))/2},{\g+\r*(sin(0)+sin(315))/2});

\draw [-{Stealth[scale=\s]}] (c2) -- ({\f+\r*cos(45)/2},{\g+\r*sin(45)/2});
\draw [-{Stealth[scale=\s]}] (c4) -- ({\f+\r*cos(135)/2},{\g+\r*sin(135)/2});
\draw [-{Stealth[scale=\s]}] (c6) -- ({\f+\r*cos(225)/2},{\g+\r*sin(225)/2});
\draw [-{Stealth[scale=\s]}] (c8) -- ({\f+\r*cos(315)/2},{\g+\r*sin(315)/2});
\draw [-{Stealth[scale=\s]}] (c0) -- ({\f+\r*cos(0)/2},{\g+\r*sin(0)/2});
\draw [-{Stealth[scale=\s]}] (c0) -- ({\f+\r*cos(90)/2},{\g+\r*sin(90)/2});
\draw [-{Stealth[scale=\s]}] (c0) -- ({\f+\r*cos(180)/2},{\g+\r*sin(180)/2});
\draw [-{Stealth[scale=\s]}] (c0) -- ({\f+\r*cos(270)/2},{\g+\r*sin(270)/2});

\end{tikzpicture}
\caption{Obstructions for outerplanar.} \label{fig:outerplanar}
\end{figure}
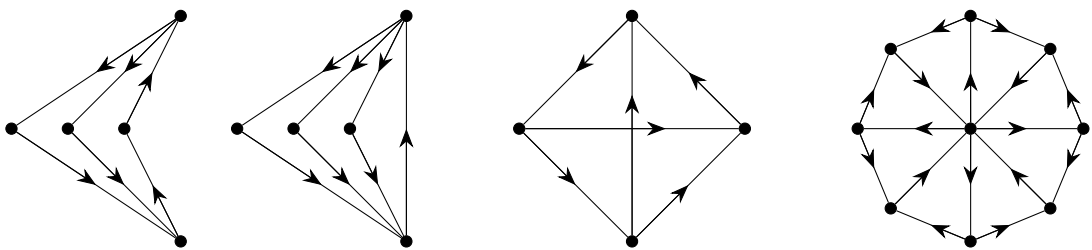

In Figure \ref{fig:outerplanar} we show four digraphs. Any directed subdivision of the first is a {\em strong theta},
and any directed subdivision of the second is a {\em reinforced theta}.
The third is called the {\em strong directed $K_4$}, and the fourth represents a class of digraphs; this one is an {\em 8-diwheel}.
In general, if $k\ge 2$ is an integer, a {\em $2k$-diwheel} is a digraph obtained from a cycle of length $2k$ by adding a new vertex adjacent with every vertex of the cycle, and then directing the edges in such a way that all the triangles become cyclic triangles.
We will show:
\begin{thm}\label{outerplanar}
A strong digraph is outerplanar if and only if no subdigraph is a directed subdivision of one of the digraphs in Figure 
\ref{fig:outerplanar}; that is, no subdigraph is a strong theta, or a reinforced theta, or a directed subdivision of the strong directed $K_4$,
 or a directed subdivision of a $2k$-wheel for some $k\ge 2$. 
\end{thm}
We remark that we could forbid subdividing the rightmost, vertical edge of the second digraph in the figure, because when this is 
subdivided, the digraph contains a strong theta. Similarly we could
forbid subdividing three edges of the strong directed $K_4$, and we could forbid subdividing the edges of the diwheel incident with the ``central'' 
vertex, for the same reason. If we do that, then we have the list of all the minimal digraphs that are strong and not outerplanar. 

For series-parallel, the list is the same except we omit the two theta-digraphs. We will show:
\begin{thm}\label{series-parallel}
A strong digraph is series-parallel if and only if no subdigraph is a directed subdivision of the strong directed $K_4$
or of a $2k$-wheel for some $k\ge 2$.
\end{thm}

The paper is organized as follows. In the next section we prove \ref{outerplanar} and \ref{series-parallel}, which are both easy. 
Then we identify all the almost-planar graphs, and this is broken into several steps. First, in section \ref{sec:4rung} we show 
that every almost-planar graph that contains a subdivision of the four-rung M\"obius ladder $V_8$ is a M\"obius chain. 
Then in section \ref{sec:8vert} we introduce two more eight-vertex graphs $U_8$ and $W_8$; and show that:
\begin{itemize}
\item the only almost-planar graph not containing a subdivision of $K_{3,3}$ is $K_5$;
\item every almost-planar graph that contains a subdivision of $K_{3,3}$ and not of $U_8,V_8$ or $W_8$ is either a M\"obius chain,
scallop or clam;
\item every almost-planar graph that contains a subdivision of $W_8$ and not of $V_8$ is a double wheel;
\item every almost-planar graph that contains a subdivision of $U_8$ and not of $V_8$ is either a double wheel, conch, whelk or mussel.
\end{itemize}
This completes the description of all almost-planar graphs. Finally, in section \ref{sec:backto} we discuss how to direct the edges of an almost-planar graph to make a minimal strong nonplanar digraph.

\section{The outerplanar and series-parallel theorems}

We will often speak of the paths and cycles of a digraph $G$. If we mean directed paths or directed cycles we will say so, and 
otherwise we mean a digraph $P$ such that $P^-$ is a path or cycle of $G^-$; it need not be directed.
The following lemma is fundamental to the paper (it is closely related to theorem 3.3 of~\cite{pfaffians}):

\begin{thm}\label{delete}
Let $G$ be a strong digraph, and let $C$ be a cycle of $G$. Then either
there is an edge $e$ of $C$ such that $G\setminus e$ 
is strong, or there is a partition $(A,B)$ of $V(G)$ such that there is only one edge from $A$ to $B$ and only one from $B$ to $A$.
\end{thm}
\Proof
For each edge $e$, let $n(e)$ be the number of directed cycles of $G$ that contain $e$, and choose $e\in E(C)$ with $n(e)$ minimum.
Let $e=ab$. If $G\setminus e$ is strong, the theorem holds, so we assume that 
$G\setminus e$ is not strong; and so there is a partition $(A,B)$ of $V(G)$ with $A,B\ne \emptyset$, such that there are no edges 
of $G\setminus e$ from $A$ to $B$. Consequently, $a\in A$ and $b\in B$, since $G$ is strong. Let $e_1\LL e_k$ be the edges of $G$
from $B$ to $A$. Every directed cycle of $G$ that contains $e_i$ also contains $e$, and contain no other edge in $\{e_1\LL e_k\}$;
and so $n(e) = n(e_1)+\cdots + n(e_k)$. Since $C$ is a cycle, one of $e_1\LL e_k$ is an edge of $C$, say $e_1$. From the choice of $e$,
$n(e_1)\ge n(e) = n(e_1)+\cdots + n(e_k)$, and so $n(e_2)\LL n(e_k)=0$. But every edge belongs to a directed cycle, since $G$ is strong,
and so $k = 1$, and the theorem holds. This proves \ref{delete}.~\bbox

\begin{thm}\label{2cut}
Let $G$ be a strong digraph, such that $G^-$ is 2-connected, and let $u,v\in V(G)$ such that $G^-\setminus \{u,v\}$ is disconnected. Let $D_1\LL D_k$
be the components of
$G^-\setminus \{u,v\}$. Then for $1\le i\le k$ there is a directed path $P_i$  of $G[V(D_i)\cup \{u,v\}]$, from $u$ to $v$ or from 
$v$ to $u$, with length at least two.
\end{thm}
\Proof There may be an edge of $G$ joining $u,v$; let $F$ be the set of all such edges. Let $G_i = G[V(D_i)\cup \{u,v\}]\setminus F$.
Since $G^-$ is 2-connected, there is a path $Q_i$ of $G_i$  between $u,v$. Let $e\in E(Q_i)$; then at least one end of $e$ is in $V(D_i)$.
But $e$ belongs to a directed cycle $C_e$ of $G$, since $G$ is strong, and if this cycle contains both $u,v$ then it contains a directed path between $u,v$
containing $e$ and the theorem holds. So we assume that for each $e\in E(Q_i)$, not both $u,v$ belong to $V(C_e)$, and so 
$C_e$ is a subdigraph of $G_i$. But then the union of the $C_e\;(e\in E(Q_i))$ is strong, and contains both $u,v$, and none of its edges 
joins $u,v$; and so it contains a directed path from $u$ to $v$ and from $v$ to $u$, both of length at least two.
This proves \ref{2cut}.~\bbox

A {\em $k$-wheel} is a graph obtained from a cycle of length $k$ by adding a new vertex adjacent to each vertex of the cycle. 
When $G$ is a graph and $F\subseteq E(G)$, an {\em $F$-cycle} of $G$ means a cycle $C$ of $G$ with $E(C)\subseteq F$.
\begin{thm}\label{4wheel}
Let $G$ be a 2-connected graph with a subgraph that is a subdivision of $K_4$. Let $F$ be the set of all edges $e$ such that 
$G\setminus e$ contains a subdivision of $K_4$ as a subgraph, and suppose there is no $F$-cycle. Then $G$ is a 
subdivision of a $k$-wheel for some $k\ge 3$.
\end{thm}
\Proof
Choose $k\ge 3$ maximum such that there is a subgraph $H$ of $G$ that is a subdivision of a $k$-wheel. Consequently
there is a cycle $C$ of $G$, and a
vertex $p_0\notin V(C)$, and
$k$ paths $P_1\LL P_k$ of $H$, each with one end $p_0$ and otherwise pairwise vertex-disjoint, and each with an end in $V(C)$ (say $p_i$)
and with no other vertex in $V(C)$.

We may assume (for a contradiction) that $H\ne G$. Since $G$ is 2-connected, there is a path $Q$ of $G$ with distinct ends $u,v\in V(H)$
and with no other vertex or edge in $H$. 
Suppose first that $u=p_0$. If $v\in V(P_i)$ for some $i$, then there is an $F$-cycle formed by $Q$ and the subpath of $P_i$ between $u,v$,
a contradiction;
so $v\in V(C)$, and $G$ contains a subdivision of a $(k+1)$-wheel, again a contradiction.
Thus $u,v\ne p_0$.
\\
\\
(1) {\em $k\ge 4$.}
\\
\\
Suppose that $k=3$. Then, from the symmetry between $p_0,p_1,p_2,p_3$, it follows that $u,v\ne p_0,p_1,p_2,p_3$. 
If there is a path $R$ of $H$ between $u,v$ containing at most one of $p_0,p_1,p_2,p_3$, 
then $E(R)\subseteq F$, 
and so $Q\cup R$ is an $F$-cycle, a contradiction. So every path of $H$ between $u,v$ passes through at least two of $p_0,p_1,p_2,p_3$.
But then adding $Q$ to $H$ makes a subdivision of $K_{3,3}$, and so all its edges are in $F$, a contradiction. This proves (1).

\bigskip

If $u,v\in V(P_1\cupcup P_k)$, let $R$ be the path between them contained in $P_1\cupcup P_k$. Since every edge of $P_1\cupcup P_k$
is in $F$, it follows that $Q\cup R$ is an $F$-cycle, a contradiction. So we assume that $u\in V(C)\setminus \{p_1\LL p_k\}$. 
Hence there is a path $R$ of $H\setminus p_0$ between $u,v$ such that
two consecutive vertices in $\{p_1\LL p_k\}$ are not in $V(R)$.
Then $E(R)\subseteq F$,
and $Q\cup R$ is an $F$-cycle, a contradiction. Thus $H=G$. This proves \ref{4wheel}.~\bbox

Consequently:
\begin{thm}\label{diwheel}
Let $G$ be a strong digraph, such that $G^-$ is 3-connected.  
Then either $G$ is a strong directed $K_4$, or $G$ is a 
$k$-diwheel for some even $k\ge 4$, or there is an edge $e$ such that $G\setminus e$ is strong and not series-parallel.
\end{thm}
\Proof
Let $F$ be the set of all edges $e$ such that
$G^-\setminus e$ contains a subdivision of $K_4$ as a subgraph.
If $C$ is a $F$-cycle,
then by \ref{delete}, $G\setminus e$ is strong for some $e\in E(C)$;
and since $G\setminus e$ is strong and not series-parallel, the third outcome holds.

So we assume that there is no $F$-cycle.
By \ref{4wheel}, $G^-$ is a subdivision of a $k$-wheel for some $k\ge 3$, and hence is a $k$-wheel, since it is 3-connected.
So 
there is a cycle $C$ of $G$ of length $k$, and a
vertex $p_0\notin V(C)$, adjacent with each vertex of $V(C)$.

If $k=3$, it follows that $G$ is a strong directed $K_4$, as required, so we assume that $k\ge 4$. 
We claim that $k$ is even and $G$ is a $k$-diwheel. Let $A$ be the set of $v\in V(G)\setminus \{p_0\}$
with zero indegree in $C$, and let $B$ be those with zero out-degree.  Thus $|A|=|B|$; if $v\in A$ then $vp_0\in E(G)$ (because $G$ is strong),
and similarly if $v\in B$ then $p_0v\in E(G)$. We need to show that $A\cup B=V(C)$. If $A=\emptyset$ (and hence $B=\emptyset$), $C$ is a
directed cycle; choose three edges of $G$ incident with $p_0$, at least one with head $p_0$ and at least one with tail $p_0$.
Then $C$ together with these three edges is strong and not series-parallel, and the third outcome holds since $k\ge 4$.
Next, suppose that
$|A|=1$ and hence $|B|=1$, and let $A=\{a\}$ and $B=\{b\}$. So $C$ is the union of two directed paths both from $a$ to $b$,
The digraph formed by the union of $C$ and the edges $bp_0, p_0a$ and a third edge incident with $p_0$ is strong and not series-parallel,
and again the third outcome holds since $k\ge 4$. So $|A|=|B|\ge 2$. The subdigraph of $G$  consisting of $C$ and all edges between 
$p_0$ and $A\cup B$ is
strong and not series-parallel, and so if $A\cup B\ne V(C)$ then the third outcome holds, and if $A\cup B=V(C)$ then $k$ is even and
$G$ is a $k$-diwheel.
This proves
\ref{diwheel}.~\bbox

Now we prove \ref{outerplanar}, which we restate:
\begin{thm}\label{outerplanar2}
A strong digraph is outerplanar if and only if no subdigraph is a 
strong theta, or a reinforced theta, or a directed subdivision of the strong directed $K_4$,
 or of a $2k$-wheel for some $k\ge 2$.
\end{thm}
\Proof Let us say a digraph is {\em forbidden} if either it is a strong theta, or a reinforced theta, or a directed subdivision of the strong directed $K_4$,
 or of a $2k$-wheel for some $k\ge 2$. The ``only if'' part of the theorem is clear, so we prove ``if''. Suppose there is a strong 
digraph $G$ that contains no forbidden digraph as a subdigraph and is not outerplanar; and choose such a digraph $G$ with 
$|V(G)|$ minimum.
\\
\\
(1) {\em $G^-$ is 3-connected.}
\\
\\
It follows easily that $G^-$ is 2-connected.
Suppose that $u,v\in V(G)$ such that $G^-\setminus \{u,v\}$ is disconnected. Let $D_1\LL D_k$
be the components of 
$G^-\setminus \{u,v\}$, and for $1\le i\le k$ let $P_i$ be a directed path of $G[V(D_i)\cup \{u,v\}]$, from $u$ to $v$ or from 
$v$ to $u$, with length at least two (this exists by \ref{2cut}). If each $P_i$ is from $u$ to $v$, then since $G$ is strong, 
there is a directed path from $v$ to $u$; and either it has length at least two, and is contained in one of the graphs 
$G[V(D_i)\cup \{u,v\}]$ (and so could replace $P_i$), or it  has length one. So we can assume that if $P_1\LL P_k$ are all from $u$ to $v$,
there is an edge $vu$; and similarly if they are all from $v$ to $u$, then there is an edge $uv$. Consequently, if $k\ge 3$
then $G$ contains a strong theta or a reinforced theta as a subdigraph, a contradiction; and so $k=2$. 

Let $P_1$ be from $u$ to $v$ say.
Let $G_1=G[V(D_i)\cup \{u,v\}]$, and let $H_1$ be obtained from $G_1$ by adding a directed edge $vu$ if it is not already present.
Thus there is a subdigraph of $G$ isomorphic to a directed subdivision of $H_1$, and so if $H_1$ contains a forbidden digraph 
then so does $G$, which is impossible. Hence $H_1$ contains no forbidden subdigraph. But $H_1$
is strong, and so is outerplanar, since $H_1$ has size less than that of $G$.  Since $D_1$ is connected, it follows that 
$H_1\setminus \{u,v\}$ is connected, and so there is an outerplanar drawing of $H_1$, such that 
the edge $vu$ is incident with the infinite region.
Similarly, choose an outerplanar
drawing of $H_2$ (defined in the same way, starting with $D_2$), with $vu$ incident with the infinite region. 
By combining these drawings, we obtain an outerplanar drawing of $G$,
a contradiction. This proves (1).

\bigskip

From \ref{diwheel},  either $G$ is a strong directed $K_4$, or $G$ is a
$2k$-diwheel for some $k\ge 2$, or there is an edge $e$ such that $G\setminus e$ is strong and not series-parallel. In the first two cases
the result holds, and in third case is impossible from the minimality of $G$.
This proves
that there is no such $G$, and so proves \ref{outerplanar2}.~\bbox

Now we prove the theorem for series-parallel digraphs, which we restate:
\begin{thm}\label{series-parallel2}
A strong digraph is series-parallel if and only if no subdigraph is a directed subdivision of the strong directed $K_4$
or of a $2k$-wheel for some $k\ge 2$.
\end{thm}
\Proof We only prove the ``if'' part, since ``only if'' is clear. 

Suppose for a contradiction that there is a strong digraph $G$ that is not series-parallel, and no subdigraph is either
a directed subdivision of either the strong directed $K_4$ or
of a $2k$-wheel for some $k\ge 2$.
Choose $G$ with $|V(G)|$  minimum. It follows that $G^-$ is 2-connected. If $G^-$ is not 3-connected, 
choose $u,v\in V(G)$ such that $G^-\setminus \{u,v\}$ is not connected. For each component $D$ of $G\setminus \{u,v\}$, as in the 
proof of \ref{outerplanar} we can add one of $uv, vu$ to $G[V(D)\cup \{u,v\}]$ to make it strong, and consequently series-parallel. 
Since this holds for each such $D$, 
it follows that $G$ is series-parallel, a contradiction. Consequently $G^-$ is 3-connected.

By \ref{diwheel},
either $G$ is a strong directed $K_4$, or $G$ is a
$k$-diwheel for some $k\ge 3$, or there is an edge $e$ such that $G\setminus e$ is strong and not series-parallel. The first 
two cases
do not hold from the choice of $G$, and
the third case is impossible from the minimality of $G$.
This proves
that there is no such $G$, and so proves 
\ref{series-parallel2}.~\bbox
\section{Almost-planar graphs containing a 4-rung M\"obius ladder}\label{sec:4rung}
Now we turn to the proof of \ref{planar}, which occupies most of the rest of the paper.
Let us say $G$ is a {\em Kuratowski digraph} if $G$ is strong and nonplanar, and minimal (under subdigraph containment) with these 
properties, and not the directed subdivision of a smaller digraph that is strong and nonplanar.
We saw in the previous section that finding the minimal digraphs that are strong and not series-parallel more-or-less reduces to find the 
3-connected graphs 
$G$ with no $F$-cycle, where $F$ is the set of edges $e$ such that $G\setminus e$ is not series-parallel. Similarly, we will see that all 
Kuratowski digraphs are orientations of almost-planar graphs.
Let us see first that:
\begin{thm}\label{reducetograph}
Let $G$ be a digraph minimal (under subdigraph) with the property that it is strong and nonplanar. Then $G$ is a directed subdivision
of a digraph $H$ such that $H^-$ is almost-planar.
\end{thm}
\Proof We proceed by induction on $|V(G)|$.
If $G$ is a directed subdivision of a digraph $H$ with $|V(H)|<|V(G)|$, then $H$ is also minimal (under subdigraph) with the property 
that it is strong and nonplanar; but the theorem holds for $H$, from the inductive hypothesis, and so also holds for $G$. Hence we
may assume that there is no such $H$, and so $G$ is a Kuratowski digraph. In this case we will prove that $G^-$ is almost-planar.

Clearly $G^-$ is 2-connected; suppose that $u,v\in V(G)$, and $G^-\setminus \{u,v\}$ is disconnected. Let $D_1\LL D_k$ be the components
of $G^-\setminus \{u,v\}$. As in the proof of \ref{outerplanar}, for $1\le i\le k$ we can add one of the edges $uv,vu$ to $G[D_i\cup \{u,v\}]$
to make it strong, and hence planar from the inductive hypothesis. But then $G$ is planar, a contradiction. This proves that
$G^-$ is 3-connected. 

Let $F$ be the set of edges $e$ such that $G\setminus e$ is 
not planar. If there is an $F$-cycle $C$ in $G$, then by \ref{delete} there is an edge $e\in E(C)$ such that $G\setminus e$ is strong; and
since $e\in F$, it follows that $G\setminus e$ is strong and nonplanar, a contradiction to the minimality of $G$. So 
$G^-$ is almost-planar. This proves \ref{reducetograph}.~\bbox

The value of this theorem is that there is something close to a converse: each of the seven types of almost-planar graphs
gives rise to infinitely many Kuratowski digraphs.  We will discuss how to 
recover a suitable orientation of $G$ later. 

In this and the next two sections we find all graphs that are almost-planar. 
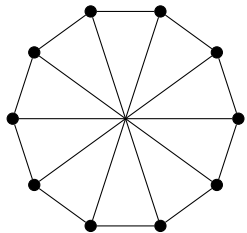
\begin{figure}[H]
\centering

\begin{tikzpicture}[scale=1.5,auto=left]
\tikzstyle{every node}=[inner sep=1.5pt, fill=black,circle,draw]
\def\r{1}
\def\a{36}
\node (a1) at ({\r*cos(0*\a)}, {\r*sin(0*\a)}){};
\node (a2) at ({\r*cos(1*\a)}, {\r*sin(1*\a)}){};
\node (a3) at ({\r*cos(2*\a)}, {\r*sin(2*\a)}){};
\node (a4) at ({\r*cos(3*\a)}, {\r*sin(3*\a)}){};
\node (a5) at ({\r*cos(4*\a)}, {\r*sin(4*\a)}){};
\node (a6) at ({\r*cos(5*\a)}, {\r*sin(5*\a)}){};
\node (a7) at ({\r*cos(6*\a)}, {\r*sin(6*\a)}){};
\node (a8) at ({\r*cos(7*\a)}, {\r*sin(7*\a)}){};
\node (a9) at ({\r*cos(8*\a)}, {\r*sin(8*\a)}){};
\node (a10) at ({\r*cos(9*\a)}, {\r*sin(9*\a)}){};

\foreach \from/\to in {a1/a6,a2/a7,a3/a8,a4/a9,a5/a10,a1/a2,a2/a3,a3/a4,a4/a5,a5/a6,a6/a7,a7/a8,a8/a9,a9/a10,a10/a1}
\draw [-] (\from) -- (\to);
\end{tikzpicture}
\caption{A five-rung M\"obius ladder.} \label{fig:Mobius}
\end{figure}

For $k\ge 2$, the {\em  M\"obius ladder} with $k$ rungs is the graph obtained from a cycle of length $2k$ by making opposite 
vertices of the cycle adjacent. Thus the 2-rung M\"obius ladder is $K_4$, the 3-rung ladder is $K_{3,3}$, and the 4-rung ladder 
is the graph known as $V_8$ or the ``Wagner graph''. Let us say the {\em ladder number} of a graph $G$ is the maximum $k$ such that 
$G$ has a subgraph that is a subdivision of the $k$-rung M\"obius ladder. (This is well-defined for almost-planar graphs since they all contain a subdivision of $K_4$.) We want to prove a theorem describing the 3-connected 
almost-planar graphs, but its proof 
breaks into cases depending on ladder number. The case
when the ladder number is two is easy, and in the case when the ladder number is at least four the argument is straightforward; but ladder number three is more tricky, and left until last.

We observe:
\begin{thm}\label{number2}
If $G$ is almost-planar, with ladder number two, then $G=K_5$.
\end{thm}
\Proof Let $G$ be almost-planar, with ladder number two. Consequently $G$ is nonplanar, and contains no subdivision of $K_{3,3}$. By a theorem of D. W. Hall~\cite{hall}, it follows that $G=K_5$. 
This proves \ref{number2}.~\bbox

\begin{figure}[H]
\centering

\begin{tikzpicture}[scale=1.5,auto=left] 
\tikzstyle{every node}=[inner sep=1.5pt, fill=black,circle,draw]
\def\r{1.5}
\def\a{18}

\node (a1) at ({\r*cos(0*\a)}, {\r*sin(0*\a)}){};
\node (a2) at ({\r*cos(1*\a)}, {\r*sin(1*\a)}){};
\node (a3) at ({\r*cos(2*\a)}, {\r*sin(2*\a)}){};
\node (a4) at ({\r*cos(3*\a)}, {\r*sin(3*\a)}){};
\node (a5) at ({\r*cos(4*\a)}, {\r*sin(4*\a)}){};
\node (a6) at ({\r*cos(5*\a)}, {\r*sin(5*\a)}){};
\node (a7) at ({\r*cos(6*\a)}, {\r*sin(6*\a)}){};
\node (a8) at ({\r*cos(7*\a)}, {\r*sin(7*\a)}){};
\node (a9) at ({\r*cos(8*\a)}, {\r*sin(8*\a)}){};
\node (a10) at ({\r*cos(9*\a)}, {\r*sin(9*\a)}){};
\node (a11) at ({\r*cos(10*\a)}, {\r*sin(10*\a)}){};
\node (a12) at ({\r*cos(11*\a)}, {\r*sin(11*\a)}){};
\node (a13) at ({\r*cos(12*\a)}, {\r*sin(12*\a)}){};
\node (a14) at ({\r*cos(13*\a)}, {\r*sin(13*\a)}){};
\node (a15) at ({\r*cos(14*\a)}, {\r*sin(14*\a)}){};
\node (a16) at ({\r*cos(15*\a)}, {\r*sin(15*\a)}){};
\node (a17) at ({\r*cos(16*\a)}, {\r*sin(16*\a)}){};
\node (a18) at ({\r*cos(17*\a)}, {\r*sin(17*\a)}){};
\node (a19) at ({\r*cos(18*\a)}, {\r*sin(18*\a)}){};
\node (a20) at ({\r*cos(19*\a)}, {\r*sin(19*\a)}){};

\foreach \from/\to in {a1/a2,a2/a3,a3/a4,a4/a5,a5/a6,a6/a7,a7/a8,a8/a9,a9/a10,
a10/a11,a11/a12,a12/a13,a13/a14,a14/a15,a15/a16,a16/a17,a17/a18,a18/a19,a19/a20, a20/a1,a1/a10, a1/a11,a1/a12,a2/a13,a3/a13,a4/a14,
a5/a14,a5/a15, a5/a16,a6/a16,a6/a17,a7/a18,a8/a18,a8/a19, a9/a20,a10/a20}
\draw [-] (\from) -- (\to);
\end{tikzpicture}
\caption{A M\"obius chain.} \label{fig:chain}
\end{figure}
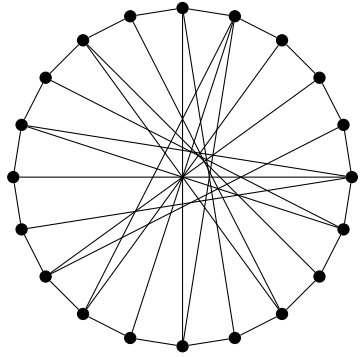

Let $C$ be a cycle of a graph $G$; then edges in $E(G)\setminus E(C)$ with both ends in $V(C)$ are called {\em chords} of $C$.
If $C$ has vertices $c_1\CC c_n\DD c_1$ in order, then
a chord with ends $c_h, c_i$ {\em $C$-crosses} one with ends $c_j, c_k$ if $h,i,j,k$ are all different and both path of $C$ 
between $c_h,c_i$ contain one of $c_j, c_k$. 
Let $G$ be a nonplanar graph obtained from a cycle $C$ by adding chords of $C$, 
subject to the rule that every two chords must either $C$-cross or share an end, and each vertex of $C$ is an end of at least one 
of the chords, and there is no cycle composed completely of chords. 
(See Figure~\ref{fig:chain}.) This graph is more general than a M\"obius ladder: let us call 
it a {\em M\"obius chain}, and $C$ is its {\em base cycle}.

For instance, the graph $K_5$
does not count as a M\"obius chain, because there is a cycle of chords. 
M\"obius chains are almost-planar. 
\begin{thm}\label{number4}
If $G$ is an almost-planar graph with ladder number at least four, then $G$ is a M\"obius chain.
\end{thm}
\Proof
Let $k$ be the ladder number of $G$, and choose a cycle $C$ of $G$ such that there are vertices $p_1\LL p_{2k}$ of $C$, distinct and
in order in $C$ (although possibly $C$ has additional vertices), such that there are paths $P_1\LL P_k$ of $G$, pairwise disjoint,
such that for $1\le i\le k$, $P_i$ has ends $p_i$ and $p_{i+k}$ and has no other vertex or edge in $C$. Let $H$ be the union of $C$ and the paths $P_1\LL P_k$; so $H$ is a subdivision of the $k$-rung M\"obius ladder. 
For $1\le i\le 2k$,
let $C_i$ be the path of  $C$ between $p_i, p_{i+1}$ that does not contain $p_{i+2}$ (reading subscripts modulo $2k$).
Let $F$ be the set of all edges $e$ such that $G\setminus e$ is nonplanar. Thus there is no $F$-cycle in $G$, since $G$ is almost-planar.
Since $k\ge 4$, $E(P_i)\subseteq F$ for $1\le i\le k$. 
Let us say an {\em $H$-jump} is a triple $(Q,u,v)$ where $u,v\in V(H)$ are distinct and $Q$ is a path of $G$ with ends $u,v$ and with no other vertex of edge in $H$. 
\\
\\
(1) {\em $u,v\in V(C)$ for every $H$-jump $(Q,u,v)$.}
\\
\\
Suppose that $u$ belongs to the interior of $P_i$ say. Then $v\notin V(P_i)$, since otherwise the subpath of $P_i$ between $u,v$ makes
an $F$-cycle with $Q$. 
There is a path $R$ of $H$ between $u,v$ consisting of a subpath of $P_i$, a subpath of $P_j$ if $v$ belongs to the interior of 
some $P_j$, and a path $R'$ of $C$ such that at most 
$(k-2)/2$  of $P_1\LL P_k$ have an end in the interior of $R'$. Consequently, the graph obtained from $H\cup Q$ by deleting the interior of $R'$
is nonplanar, and so $E(R')\subseteq F$. It follows that $E(R)\subseteq F$, and therefore $Q\cup R$ is an $F$-cycle, a contradiction.
This proves (1).
\\
\\
(2) {\em For every $H$-jump $(Q,u,v)$ and every path $R$ of $C$ between $u,v$, at least $k-2$ of $P_1\LL P_k$ have an end in the interior of $R$.}
\\
\\
If not, then the graph obtained from $H\cup Q$ by deleting the interior of $R$ is nonplanar, and so $E(R)\subseteq F$, and $Q\cup R$ is an $F$-cycle, which is impossible. This proves (2).
\\
\\
(3) {\em For every $H$-jump $(Q,u,v)$, every path $R$ of $C$ between $u,v$ contains an end of each of $P_1\LL P_k$.}
\\
\\
Suppose not; then from (2), not both $u,v$ are ends of $P_1\LL P_k$, so we assume that $u\notin \{p_1\LL p_{2k}\}$. 
We assume that $u$
belongs to the interior of $C_{2k}$, and $p_1\LL p_{k-2}$ belong to the interior of $R$. The vertex $p_{k-1}$ might belong 
to the interior of $R$, or be an end of $R$ (and hence equal $v$), or not belong to $R$. But certainly $p_1,p_2$ belong to the 
interior of $R$. The graph obtained from $H\cup Q$ by deleting the interior of $C_1$ is nonplanar; indeed, if $j\le k-1$ is 
maximum such that $p_j\in V(R)$, then the graph obtained from $H\cup Q$ by deleting the interior of $C_1\cup C_2\cupcup C_{j-1}$ 
is nonplanar. Consequently $E(C_1)\subseteq F$. But the graph obtained from $H\cup Q$ by deleting the interior of $C_{k+1}$
is topologically equivalent to the one obtained by deleting the interior of $C_1$, and so is also nonplanar; and so 
$E(C_{k+1})\subseteq F$. Then $C_1\cup C_{k+1}\cup P_1\cup P_2$  is an $F$-cycle, a contradiction. This proves (3).
\\
\\
(4) {\em For every $H$-jump $(Q,u,v)$, exactly one of $u,v$ is in $\{p_1\LL p_{2k}\}$.}
\\
\\
If neither of $u,v$ belong to $\{p_1\LL p_{2k}\}$, then from (3), $H\cup Q$ is a subdivision of a $(k+1)$-rung M\"obius ladder, 
contradicting that $k$ is the ladder number of $G$. If both $u,v\in \{p_1\LL p_{2k}\}$, then by (3) $u,v$ are ends of the same one of $P_1\LL P_k$,
which contradicts that there is no $F$-cycle. This proves (4).
\\ 
\\
(5) {\em Each of $P_1\LL P_k$ has length one, and $V(H) = V(G)$.}
\\
\\
Suppose that there is a component $D$ of $G\setminus V(H)$. Since $G$ is 3-connected, there are three vertices in $V(H)$ that
have a neighbour in $V(D)$, and every two of them can be joined by a path with the properties of the path $Q$. In each case
one of these three paths violates (4). So $V(H)=V(G)$. If some $P_i$ has length more than one, then since $G$ is 3-connected,
there is an $H$-jump $(Q,u,v)$ such that one of $u,v$ is in the interior of $P_i$, contrary to (1). This proves (5).

\bigskip

There is no cycle with no edges in $E(C)$, since any such cycle would be an $F$-cycle, a contradiction.
From (5), $G$ is obtained from $H$ by adding chords $C$; and by (3), every chord $C$-crosses or meets each of 
the edges $p_1p_{k+1}\LL p_kp_{2k}$. To complete the proof, we must show that every two additional chords meet or $C$-cross each other. 
Thus, let $p_{2k}a$ and $p_ib$ be edges of $G$ not in $E(H)$. (This is without loss of generality, since each such edge has an end in 
$\{p_1\LL p_{2k}\}$ by (4).) From (5), $a,b\notin \{p_1\LL p_{2k}\}$. 
From (3), $a\in V(C_{k-1}\cup C_k)$ and from (5), $a\notin  \{p_1\LL p_{2k}\}$, so $a$ belongs to the interior of one of $C_{k-1}, C_k$. Let $H'$ be obtained from $H$ by deleting the edge $p_kp_{2k}$ and adding $p_{2k}a$; this is another $k$-rung M\"obius ladder 
with the same base cycle $C$. From (3) applied to $H'$, $p_ib$ meets or $C$-crosses $p_{2k}a$. 
This proves \ref{number4}.~\bbox

\section{Excluding subdivisions of $U_8$, $V_8$ and $W_8$}\label{sec:8vert}
Next we turn to the case of almost-planar graphs with ladder number three, that is, when $G$ contains a subdivision of $K_{3,3}$ but not of $V_8$. It turns out that there are several kinds of almost-planar graph with this property, and in this section we look at those
that also contain no subdivision of two other graphs $U_8,W_8$, defined later.

We shall need to look at the structure of $G$ relative 
to a subdivision of $K_{3,3}$ it contains; let us set up some notation for this. If $H$ is a subgraph of $G$ that 
is a subdivision of $K_{3,3}$, there are distinct $a_1,a_2,a_3,b_1,b_2,b_3\in V(G)$ and paths $P_{i,j}$ between $a_i$ and $b_j$ for all $i,j\in \{1,2,3\}$, 
pairwise vertex-disjoint except for their ends, such that $H$ is the union of these nine paths. Let us call this the {\em standard notation} for $H$.

\begin{thm}\label{number31}
Let $G$ be an almost-planar graph with ladder number three, and let $H$ be a subgraph that is a subdivision of $K_{3,3}$.
Then $V(H)=V(G)$. Moreover, if $e=uv$ is an edge of $G\setminus E(H)$, then (in the standard notation)  $e$ is incident with at least 
one of 
$\{a_1,a_2,a_3,b_1,b_2,b_3\}$, and there do not exist $i,j\in \{1,2,3\}$ such that $u,v\in V(P_{i,j})$.
\end{thm}
\Proof
Let $F$ be the set of all edges $e$ such that $G\setminus e$ is nonplanar. We use the standard notation for $H$.
\\
\\
(1) {\em If $Q$ is a path between 
distinct $u,v\in V(H)$, with no other vertex or edge in $H$, then 
there do not exist $i,j\in \{1,2,3\}$ such that $u,v\in V(P_{i,j})$; and 
at least one of $u,v$ is in $\{a_1,a_2,a_3,b_1,b_2,b_3\}$.}
\\
\\
If $u,v\in V(P_{i,j})$ then the subpath of $P_{i,j}$
joining $u,v$, together with $Q$, is an $F$-cycle, so there is no such $(i,j)$.
For the second statement, we may assume that $u$ belongs to the interior of $P_{1,1}$. 
and $v$ belongs to the interior of $P_{i,j}$
for some $(i,j)\ne (1,1)$. If $i=1$, let $R$ be the path of $H$ between $u,v$ consisting of subpaths of $P_{1,1}$ from $u$ to $a_1$
and a subpath of $P_{1,j}$ from $a_1$ to $v$; then all its edges are in $F$, a contradiction. So $i\ne 1$ and similarly $j\ne 1$.
But then $H\cup Q$ is a subdivision of a four-rung M\"obius ladder, a contradiction. This proves (1).

\bigskip
Now suppose that $V(H)\ne V(G)$, and let $D$ be a component of $G\setminus V(H)$. Since $G$ is 3-connected, there are at least three 
vertices in $V(H)$ that have a neighbour in $V(D)$, say $u,v,w$. By the second statement of (1), at least two of $u,v,w$ belongs to $\{a_1,a_2,a_3,b_1,b_2,b_3\}$;
and by the first statement of (1), one of $\{a_1,a_2,a_3\}, \{b_1,b_2,b_3\}$ contains none of $u,v,w$. Thus, from the symmetry and (1), 
we may assume that $v=a_2, w=a_3$, and
either $u=a_1$ or $u$ belongs to the interior of $P_{1,1}$. In either case, adding $D$ and the edges between $D,H$ to $H$, and 
then deleting $b_2$ and the interiors of $P_{1,2}, P_{2,2}, P_{3,2}$ makes a nonplanar graph; so all edges of $P_{1,2}, P_{2,2}, P_{3,2}$ belong to $F$.
Similarly all edges of $P_{1,3}, P_{2,3}, P_{3,3}$ belong to $F$, and so there is an $F$-cycle, a contradiction. Thus $V(H)=V(G)$. This proves \ref{number31}.~\bbox

We need to define three further types of graph. 
Let us say a {\em scallop} is a 3-connected graph either equal to $K_5$, or consisting of a cycle with vertices $c_1\CC c_n\DD c_1$ in order (with
$n\ge 3$), and three further vertices 
$u,v,w$; $u,v,c_1,c_n$ are pairwise adjacent, and $w$ is adjacent to $u,v$ and to $c_2\LL c_{n-1}$, and possibly to $c_1,c_n$. 

A {\em clam} is a 3-connected graph consisting of a cycle $C$ with vertices $c_1\CC c_n\DD c_1$ in order,
and two further vertices $u,v$, adjacent, with the following properties:
\begin{itemize}
\item there exists $j\in \{2\LL n-1\}$ such that $u,v$ are each
adjacent to each of $c_1,c_j,c_n$;
\item the set of neighbours of $u$ in $V(C)$ is $\{c_n,c_1,c_2\LL c_j\}$, and the set of neighbours of $v$ in $V(C)$ is
$\{c_j, c_{j+1}\LL c_n,c_1\}$.
\end{itemize}
A {\em whelk} is similar. There is a cycle $C=c_1\CC c_n\DD c_1$ and two extra vertices $u,v$, adjacent to each other.
There exist $i,j$ with $2\le i<j\le n$ such that the set of neighbours of $u$ in $V(C)$ is $\{c_1,c_2\LL c_j\}$, and the set of neighbours of $v$ in $V(C)$ is
$\{c_j, c_{j+1}\LL c_n,c_i\}$.
(See Figure \ref{fig:clam}.) Clams and whelks are both almost-planar.

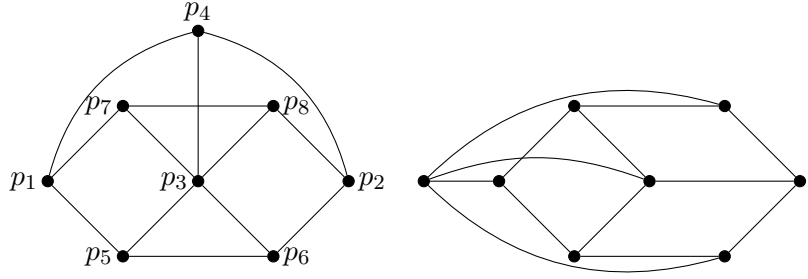
\begin{figure}[H]
\centering

\begin{tikzpicture}[scale=1,auto=left]
\tikzstyle{every node}=[inner sep=1.5pt, fill=black,circle,draw]
\node (a1) at (-2,0){};
\node (a2) at (2,0){};
\node (a3) at (0,0){};

\node (a4) at (-1,-1){};
\node (a5) at (1,-1){};
\node (a6) at (-1,1){};
\node (a7) at (1,1){};
\node (a8) at (0,2){};

\node (b1) at (3,0){};
\node (b2) at (4,0){};
\node (b3) at (5,1){};
\node (b4) at (5,-1){};
\node (b5) at (6,0){};
\node (b6) at (7,1){};
\node (b7) at (7,-1){};
\node (b8) at (8,0){};
\tikzstyle{every node}=[]

\draw (a1) node [left]           {$p_1$};
\draw (a2) node [right]           {$p_2$};
\draw (a3) node [left]           {$p_3$};
\draw (a4) node [left]           {$p_5$};
\draw (a5) node [right]           {$p_6$};
\draw (a6) node [left]           {$p_7$};
\draw (a7) node [right]           {$p_8$};
\draw (a8) node [above]           {$p_4$};

\foreach \from/\to in {a1/a6,a1/a4,a2/a5,a2/a7,a3/a4,a3/a5,a3/a6,a3/a7, a3/a8,a4/a5,a6/a7, b1/b2,b2/b3,b2/b4,b3/b5,b4/b5,b3/b6,b4/b7,b5/b8,b6/b8,b7/b8}
\draw [-] (\from) -- (\to);
\draw (a1) to [bend left = 30] (a8);
\draw (a2) to [bend right = 30] (a8);
\draw (b1) to [bend left = 30] (b6);
\draw (b1) to [bend right = 30] (b7);
\draw (b1) to [bend left = 20] (b5);

\end{tikzpicture}
\caption{The graph $U_8$ and $W_8$.} \label{fig:anti}
\end{figure}

The graphs $U_8, W_8$ are shown in Figure \ref{fig:anti}.
They are both 
almost-planar, with ladder number three, and not M\"obius chains, and they are the only minimal such graphs:
\begin{thm}\label{noU8}
Let $G$ be an almost-planar graph with ladder number three that contains no subdivision of $U_8$ or $W_8$. Then 
$G$ is either a M\"obius chain or a scallop or clam.
\end{thm}
\Proof
Choose a subgraph $H$ of $G$ that is a subdivision of $K_{3,3}$, using the standard notation. By \ref{number31}, $V(H)=V(G)$, and 
every edge of $G$ not in $E(H)$ has at least one end in $\{a_1,a_2,a_3,b_1,b_2,b_3\}$. We need to look at such edges with {\em exactly} 
one end in this set; we will worry about those with both ends in the set later. Let us say that $(a_i,b_j)$ is {\em active}
if $a_i$ is adjacent to a vertex in the interior of $P_{i',j}$ for some $i'\in \{1,2,3\}$ different from $i$. We define that 
$(b_j,a_i)$ is active similarly; that is, if $b_j$ is adjacent to a vertex in the interior of $P_{i,j'}$ for some $j'\in \{1,2,3\}$ 
different from $j$. 
\\
\\
(1) {\em If $(a_i,b_j)$ is active then $P_{i,j}$ has length one.}
\\
\\
Without loss of generality, let $a_1$ have a neighbour $v$ in the interior of $P_{2,1}$. Let $H'$ be be the subdivision of $K_{3,3}$ 
obtained from $H$ by adding the edge $a_1v$ and deleting the interior of $P_{1,1}$; then $V(H') = V(G)$ by \ref{number31},
and so $P_{1,1}$ has length one. This proves (1).
\\
\\
(2) {\em For all $i,j\in \{1,2,3\}$, not both $(a_i,b_j)$ and $(b_j,a_i)$ are active.}
\\
\\
Suppose (without loss of generality) that $a_1$ has a neighbour in the interior of $P_{2,1}$, and $b_1$ has a neighbour in the interior of $P_{1,2}$.
The graph obtained from $H$ by adding these two edges and deleting the interior of $P_{1,1}$ is a subdivision of $V_8$, a contradiction. 
This proves (2).
\\
\\
(3) {\em If $i,j,j'\in \{1,2,3\}$ with $j\ne j'$, then not both $(a_i,b_j)$ and $(a_i, b_{j'})$ are active. 
The same holds with the $a_i$'s and $b_j$'s exchanged.}
\\
\\
If both pairs are active, the graph obtained from $H$ by adding the corresponding two edges is a subdivision of $U_8$, 
a contradiction. This proves (3).
\\
\\
(4) {\em If $i,j,j'\in \{1,2,3\}$ with $j\ne j'$, then not both $(a_i,b_j)$ and $(b_{j'}, a_i)$ are active. 
The same holds with the $a_i$'s and $b_j$'s exchanged.}
\\
\\
Suppose (without loss of generality) that $a_1$ has a neighbour $u$ in the interior of $P_{2,1}$, and $b_2$ has a neighbour $v$ in the interior
of one of $P_{1,1}, P_{1,3}$. By (1), 
$P_{1,1}$ has length one, and so $v$ is in the interior of $P_{1,3}$. Let $H'$
be the subdivision of $K_{3,3}$ obtained from $H$ by adding the edge $b_2v$ and deleting the interior of $P_{1,2}$; then the edge $a_1u$
violates \ref{number31} applied to $H'$. This proves (4).
\\
\\
(5) {\em If $i,j,j'\in \{1,2,3\}$ with $j\ne j'$, then not both $(b_j, a_i)$ and $(b_{j'}, a_i)$ are active. 
The same holds with the $a_i$'s and $b_j$'s exchanged.}
\\
\\
Suppose (without loss of generality) that $(b_1,a_1)$ and $(b_2,a_1)$ are active. By (1), $P_{1,1}$ and $P_{1,2}$ both have length one,
and so $b_1, b_2$ both have a neighour in the interior of $P_{1,3}$. Let $u,v$ be neighbours of $b_1,b_2$ respectively in the 
interior of $P_{1,3}$. By exchanging $b_1,b_2$ if necessary, we may assume that $u$ lies in the subpath of $P_{1,3}$
between $a_1$ and $v$. But then the edge $b_1u$ violates \ref{number31} applied to the subdivision of $K_{3,3}$ obtained from 
$H$ by adding the edge $b_2v$ and deleting the edge $a_1b_2$.  This proves (5).

\bigskip
From (2)--(5), there are at most three active pairs, and they are vertex-disjoint in the natural sense. Thus we may assume that 
there are no active pairs except possibly one of $(a_1,b_1), (b_1,a_1)$, one of $(a_2,b_2),(b_2,a_2)$ and one of $(a_3,b_3),(b_3,a_3)$.
\\
\\
(6) {\em $P_{1,1}, P_{2,2}, P_{3,3}$ each have length one.}
\\
\\
If the interior of $P_{1,1}$ is nonempty, then since $G$ is 3-connected, there is a path $Q$ of $G$ with one end in this interior,
and the other in $V(H)\setminus V(P_{1,1})$, and with no other vertices or edges in $H$. By \ref{number31}, $Q$ has length one and its
end outside $P_{1,1}$ is one of $a_2,a_3,b_2,b_3$, contradicting that no such pair is active. This proves (6).

\bigskip

Let $C$ be the cycle 
$$P_{1,3}\cup P_{2,3}\cup P_{2,1}\cup P_{3,1}\cup P_{3,2}\cup P_{1,2}.$$
Thus, $V(C)=V(G)$, from (6) and since $V(H)=V(G)$.
Let $\mathcal{A}$ be the set of edges of $G$ with both ends in $\{a_1,a_2,a_3\}$, let $\mathcal{B}$ be the set with both
ends in $\{b_1,b_2,b_3\}$, and let $\mathcal{E}$ be the set of edges in $E(G)\setminus E(H)$ that are not in 
$\mathcal{A}\cup \mathcal{B}$. 
\\
\\
(7) {\em Every two edges in $\mathcal{E}$ either meet or $C$-cross.}
\\
\\
Let $a_1$ have a neighbour $u$ in the interior of $P_{2,1}$; it suffices to show that $a_1u$ meets or $C$-crosses every other edge
in $\mathcal{E}$. 
Certainly $a_1u$ meets or crosses every edge in $\mathcal{E}$ incident with
$a_1, a_2, a_3$ or $b_3$, and so it remains to handle those
incident with $b_1$ and those incident with $b_2$. By (2), there are no edges of $\mathcal{E}$ incident with $b_1$. 
Suppose that $b_2v\in \mathcal{E}$ does not meet or cross $a_1u$.
Then $v$ 
belongs to the subpath of $P_{2,1}$ between $b_1$ and $u$, and $v\ne u$. But then 
the graph obtained from $H$ by adding the edges $a_1u$ and $b_2v$ to $H$ is a subdivision of $W_8$, which is impossible.
This proves (7).

\bigskip

So far, these results apply to any choice of $H$, but now let us choose $H$ to maximize the number of active pairs. Since now
we will need to compare the number of active pairs for different choices of $H$, from now on we will write ``$H$-active'' for active.
\\
\\
(8) {\em Every edge in $\mathcal{E}$ meets or $C$-crosses
every edge in $\mathcal{A}\cup \mathcal{B}$.}
\\
\\
Let $a_1$ have a neighbour $u$ in the interior of $P_{2,1}$; it suffices to show that $a_1u$ meets or $C$-crosses every edge
in $\mathcal{A}\cup \mathcal{B}$. This is clear for all such edges except possibly $b_1b_2$, so we assume that $b_1b_2\in E(G)$.
Let $H'$ be obtained from $H$ by deleting $a_1b_1$ and adding $a_1u$.
Then $(b_2,a_3)$ is $H'$-active (because of $b_1b_2$), and $(a_1,u)$ is $H'$-active (because of $a_1b_1$).
The pair $(b_1,a_1)$ is not $H$-active, by (2);
if either of $(a_3,b_3),(b_3,a_3)$ is $H$-active then it is also $H'$-active; if $(a_2,b_2)$ is $H$-active then it is $H'$-active;
so from the choice of $H$ (maximizing the number of active pairs), it follows that $(b_2,a_2)$ is $H$-active and not $H'$-active.
Hence there is a neighbour $v$ of $b_2$ that belongs to the subpath of $P_{2,1}$ between $u,b_1$, with $v\ne b_1$. 
Since $b_2v\in \mathcal{E}$, from (7) it follows that $v=u$. But then $a_1\DD u\DD b_2\DD b_1\DD a_1$ is an $F$-cycle, 
a contradiction. This proves (8).
\\
\\
(9) {\em If some two edges in $\mathcal{A}\cup \mathcal{B}$ do not meet or $C$-cross, then $G$ is a scallop.}
\\
\\
Certanly every two edges in $\mathcal{A}$ meet, and the same for $\mathcal{B}$; so we may assume that $a_1a_2,b_1b_2\in E(G)$. 
Suppose that $(a_1,b_1)$ is $H$-active, and let $a_1$ have a neighbour $u$ in the interior of one of $P_{2,1}, P_{3,1}$. By (8),
$u\notin V(P_{2,1})$ and so $u$ is in the interior of $P_{3,1}$. But then $a_1\DD b_1\DD b_2\DD a_2\DD a_1$ is an $F$-cycle,
a contradiction. So $(a_1,b_1)$ is not $H$-active, and similarly none of $(b_1.a_1), (a_2,b_2), (b_2,a_2)$ are $H$-active. By (2) and the symmetry,
we may assume
that $(a_3,b_3)$ is not $H$-active, so all edges in $\mathcal{E}$ are incident with $b_3$. If neither of $a_1a_3,a_2a_3$ exists, then $G$ is a
scallop,
so we assume that $a_1a_3$ exists. If $\mathcal{E}\ne \emptyset$, let $b_1v\in \mathcal{E}$; so $v$ belongs to the interior of one of $P_{3,1},
P_{3,2}$. By (8), since $b_1v$ meets or $C$-crosses $a_1a_3$, it follows that $v\in V(P_{3,2})$; but then $a_1\DD a_2\DD b_2\DD b_1\DD a_1$ is an 
$F$-cycle, a contradiction. So $\mathcal{E}=\emptyset$, and so $|V(G)|=6$. If neither of the edges $b_1b_3,b_2b_3$ exists, then $G$ is a scallop; 
and since there is symmetry exchanging $b_1,b_2$, we may assume that $b_1b_3\in E(G)$. But then every edge between $\{a_1,b_1\}$ and
$\{a_2,a_3,b_2,b_3\}$ is in $F$, and so there is an $F$-cycle, a contradiction. This proves (9).

\bigskip
In view of (7), (8) and (9), we assume that every two edges in $E(G)\setminus E(C)$ meet or $C$-cross. 
If there is no cycle of chords for $C$, then $G$ is a M\"obius chain with base cycle $C$; so we assume that there is a cycle $D$ edge-dsjoint from
$C$. In this case we will prove that $G$ is a clam.
Since $D$ is not an $F$-cycle, there is an edge $f\in E(D)\setminus F$. 
Since $f\notin F$, it follows that $f=a_ib_i$ for some $i\in \{1,2,3\}$,
so we may assume that $f=a_1b_1$. Since $f\notin F$, $(a_1,b_1)$ and $(b_1,a_1)$ are not $H$-active; and since $f\in E(D)$,
we may assume that $a_1a_3, b_1b_2\in E(D)$. Let $J$ be the graph obtained from $H$ by adding
$a_1a_3, b_1b_2$, and let us examine the edges in $G\setminus E(J)$.
\begin{itemize}
\item If $(a_2,b_2)$ is $H$-active, $a_2$ has a neighbour in the interior of $P_{3,2}$ (it cannot be in the interior of $P_{1,2}$ 
since it must $C$-cross $b_1b_2$). But then $G\setminus f$ is nonplanar, a contradiction. So $(a_2,b_2)$ is not $H$-active. Since 
$a_2a_1\notin E(G)$ (because it does not $C$-cross $b_1b_2$), it follows that the only possible edge in $E(G)\setminus E(J)$ incident with $a_2$
is $a_2a_3$. 
Similarly no edge of $G\setminus E(J)$ is incident with $b_3$ except possible $b_2b_3$.
\item $(a_1,b_1)$ and $(b_1,a_1)$ are not $H$-active, since $a_1b_1\notin F$. So there are no edges of $G\setminus E(J)$ incident with 
$a_1$ or $b_1$. 
\end{itemize}
Consequently, every edge of $G\setminus E(J)$ is incident with one of $a_3,b_2$, and so 
$D$ has length five, and there exists $c\in V(P_{2,3})$ such that $D$ is the cycle
$a_1\DD a_3\DD c\DD b_2\DD b_1\DD a_1$. Since all edges of $G\setminus E(C)$ incident with $b_2$ meet or $C$-cross both 
$a_1b_1$ and $a_3c$,
it follows that all remaining neighbours of $b_2$ belong to $P_{2,1}$ or to the subpath of $P_{2,3}$ between $a_2$ and $c$, and
a similar statement holds for $a_2$. Consequently $G$ is a clam.
This proves \ref{noU8}.~\bbox
\section{Double wheels}\label{sec:doublewheels}

That concludes the ``M\"obius''-type results. Now we look at the almost-planar graphs with ladder number three that contain a subdivision of one of $U_8, W_8$.
Let us say a graph $G$ is a {\em double wheel} if it is 3-connected and nonplanar, and consists of a cycle $C$ and two further vertices $u,v$, adjacent to each 
other, such that  
at most one vertex of $C$ is adjacent to both $u,v$. 
(Thus clams and whelks are not double wheels.)
We call $u,v$ its {\em wheel centres} and $C$ its {\em base cycle}.
Double wheels are almost-planar.

\begin{thm}\label{doublewheel}
Let $G$ be an almost-planar graph with ladder number three, that contains a subdivision of $W_8$. Then $G$ is a double wheel.
\end{thm}
\Proof
Let $C$ be a cycle of $G$, and let $a,b\in V(G)\setminus V(C)$; let $p_1,p_2,p_3,p_4,p_5,p_6\in V(C)$ be in cyclic order;
and let $P_1\LL P_6$ paths between $\{a,b\}$ and $V(C)$, with no internal vertex in $\{a,b\}\cup V(C)$, and pairwise vertex-disjoint except for their ends, where $P_i$ is between $a$ and $p_i$ for $i = 1,3,5$ and $P_i$ is between $b$ and $p_i$ for $i = 2,4,6$. Let $P_0$
be a path between $a,b$, vertex-disjoint from $C$ and with no internal vertex in $P_1\cupcup P_6$. For $1\le i\le 6$ let $C_i$
be the path of $C$ between $p_i$ and $p_{i+1}$ (reading subscripts modulo 6). Let $H$ be the union of $C$ and $P_0,P_1\LL P_6$.
Since the graph obtained from $H$ by deleting the interior of $P_i$ contains a subdivision of $K_{3,3}$, \ref{number31} implies
that $P_i$ has length one, for $1\le i\le 6$, and $V(H)=V(G)$. As usual, let $F$ be the set of edges $e$ such that $G\setminus e$
is nonplanar; so $ap_1,ap_3,ap_5,bp_2,bp_4,bp_6\in F$. 

If every edge of $G\setminus E(H)$ is incident with $a$ or $b$, then $P_0$ has length one (since $G$ is 3-connected) and so $G$ is a double wheel.
Hence we suppose (for a contradiction), that 
$uv\in E(G)\setminus E(H)$, and $u,v\ne a,b$. 
\\
\\
(1) {\em Every edge of $G\setminus E(H)$ is incident with $a$ or $b$, and $P_0$ has length one.}
\\
\\
Suppose first that $P_0$ has legth more than one.
Then there is an edge $uv\in E(G)\setminus E(H)$
where $u$ belongs to the interior of $P_0$. Thus $v\notin V(P_0)$ (since otherwise there would be an $F$-cycle,
as usual), so $v\in V(C)$. From the symmetry we may assume that either $v=p_1$ or $v$ belongs to the interior of $C_1$. In either case,
adding $uv$ to $H$ and deleting $ap_1$ and $bp_2$ gives a subdivision of $V_8$, a contradiction. This proves that  
$P_0$ has length one. Now suppose that there is an edge $uv\in E(G)\setminus E(H)$, and $u,v\ne a,b$. It follows that $u,v\in V(C)$.
From \ref{number31}, one of $u,v$ equals one of $p_1\LL p_6$, say $u = p_1$ without loss of generality. By \ref{number31}
again, applied to the subdivision of $K_{3,3}$ obtained from $H$ by deleting $ap_1$ and $bp_6$, $v$ is one of $p_2,p_3,p_4,p_5$.
Since there is no $F$-cycle, $v\ne p_2$.If $v=p_4$ then $ab\in F$, and $a\DD p_1\DD p_4\DD b\DD a$ is an $F$-cyce, a contradiction.
So 
$v\in \{p_3,p_5\}$ and we may assume that $v=p_3$ from the symmetry. But the graph obtained from
$H$ by adding $p_1p_3$ and deleting $p_2$ is nonplanar, and so $E(C_1\cup C_2)\subseteq F$, and since $p_1p_3\in F$
this is impossible. This proves (1).

\bigskip

From (1), to see that $G$ is a double wheel, note that every edge between $\{a,b\}$ and $V(C)$
belongs to $F$, and so $a,b$ have at most one common neighbour in $C$ since there is no $F$-cycle. This proves \ref{doublewheel}.~\bbox


A {\em conch} is a 3-connected graph $G$ with the following properties. There are three vertices $p_1,p_2,p_3$, and three paths $P,Q,R$,
vertex-disjoint except for their ends. The paths $P, Q$ are both between $p_1,p_2$, and do not contain $p_3$, and $R$ is between 
$p_2,p_3$, and does not contain $p_1$. The vertex $p_3$
has at least two neighbours in the interior of each of $P,Q$, and $p_1$ has at least two neighbours in the interior of $R$.
There are no other edges, except possibly $p_1$ is adjacent to $p_2$ or $p_3$. 

A {\em mussel} is a 3-connected graph with the following properties. There are three vertices $p_1,p_2,p_3$ and three paths $P,Q,R$
between $p_1,p_2$, pairwise vertex-disjoint except for their ends, and each of length at least three. The vertex $p_3$ belongs to none of these paths, but is adjacent to every vertex in their interiors, 
and may also be adjacent to $p_1$ or to $p_2$. Moreover, $p_1,p_2$ may be adjacent, although at most two of the edges $p_1p_2,p_2p_3,p_1p_3$ are present. 
Let us call $p_1,p_2$ its {\em tips} and $p_3$ its {\em hinge}.
(See Figure \ref{fig:clam}.) It is easy to check that all conches and mussels are almost-planar.
\begin{thm}\label{conch}
Let $G$ be an almost-planar graph with ladder number three, that contains a subdivision of $U_8$. Then $G$ is 
a double wheel, conch, whelk, or mussel.
\end{thm}
\Proof
Since $G$ contains a subdivision $H$ of $U_8$, we may choose vertices $p_1\LL p_8$, and paths $P_{i,j}$ joining $p_i,p_j$ for each edge $p_ip_j$
of $U_8$, labelled as in Figure \ref{fig:anti}, where $H$ is the union of all these paths. There is a subdivision of $K_{3,3}$ 
obtained by deleting from $H$ the interiors of $P_{3,5}$ and $P_{3,7}$, and so by \ref{number31}, $V(H)=V(G)$, and 
$P_{3,5}$ and $P_{3,7}$ have 
length one, and similarly $P_{3,6}$ and $P_{3,8}$ have length one. Let $F$ be the set of edges $e$ such that $G\setminus e$ is nonplanar.
Thus $p_3p_5,p_3p_6,p_3p_7,p_3p_8\in F$. 
\\
\\
(1) {\em For each edge $uv\in E(G)\setminus E(H)$, one of $u,v\in \{p_1,p_2,p_3,p_4\}$.}
\\
\\
By \ref{number31} applied to the same subdivision of $K_{3,3}$, it follows that one of $u,v\in \{p_1,p_2,p_3,p_4,p_6,p_8\}$, and so 
we may assume that $u=p_6$ without loss of generality. From \ref{number31} applied to the subdivision of $K_{3,3}$ 
obtained from $H$ by deleting the edges $p_3p_6$ and $p_3p_7$ , it follows that $v\in \{p_1,p_2,p_3,p_4,p_5,p_8\}$, so we may assume that
$v\in \{p_5,p_8\}$; and $v\ne p_8$ by \ref{number31} applied to the subdivision of $K_{3,3}$
obtained by deleting from $H$ the interiors of $P_{3,6}$ and $P_{3,8}$, Hence $v=p_5$, contradicting that there is no $F$-cycle. This 
proves (1).
\\
\\
(2) {\em If no edge of $G\setminus E(H)$ is incident with $p_1$ or with $p_2$ then $G$ is either a double wheel, a whelk, 
 or a mussel.}
\\
\\
In this case, $P_{3,4}$ has length one, since no edge of $G\setminus E(H)$ has an end in its interior. Suppose first that 
$p_3$ has no neighbours
in the interior of $P_{1,4}$ or $ P_{2,4}$. We claim that either $G$ is a double wheel, or a whelk.
To show this, we assume it is not a double wheel, and so there are at least two vertices in $V(C)$ adjacent to both $p_3,p_4$, say $c,d$. 
All edges between $\{p_3,p_4\}$ and $V(C)$ belong to $F$ except possibly $p_1p_4,p_2p_4$; so, since there is no $F$-cycle, 
we may assume that $d=p_1$  and $P_{1,4}$ has length one, and $p_1p_4\notin F$.
Since $G\setminus p_1p_4$ is planar, there is a subpath $Q$ of $P_{2,8}\cup P_{2,6}$ containing all neighbours of $p_4$ in $V(C)$
except $p_1$, and $p_3$ has no neighbour in the interior of $Q$. Thus $c$ is one end of $Q$, and hence $c\in V(P_{2,8}\cup P_{2,6})$. 
If there is a third vertex $e$ say
of $C$ that is adjacent to both $p_3,p_4$, then $e$ is necessarily the other end of $Q$, and so also belongs to $V(P_{2,8}\cup P_{2,6})$.
By the same argument, applied to $c,e$, it follows that one of $c,e$ equals $p_1$, a contradiction. So there is no such $e$, and so
$G$ is a whelk.

Now we assume that $p_3$ has a neighbour $u$ in the interior of 
one of $P_{1,4}$, $ P_{2,4}$. If $p_4$ is incident with no edge of $G\setminus E(H)$ then $G$ is a mussel, so we assume that
$p_4v\in E(G)\setminus E(H)$. From the symmetry we may assume that $v$ is in the interior of the path $P_{1,5}\cup P_{5,6}$. But then adding the edge $p_3u$ to $H$ and deleting $p_3p_4, p_3p_5$ and $p_3p_7$ gives a subdivision of $K_{3,3}$ in which the edge $p_4v$
violates \ref{number31}. This proves (2).

\bigskip

From (2) we may assume that there is an edge $p_1u\in E(G)\setminus E(H)$.
\\
\\
(3) {\em $u\in V(P_{3,4}\cup P_{2,4})$, and $P_{1,4}$ has length one. }
\\
\\
If $v$ belongs to $P_{1,4}\cup P_{1,5}\cup P_{5,6}\cup P_{1,7}\cup P_{7,8}$, there is an $F$-cycle using $e$, a contradiction.
So we may assume that $v$ belongs to the interior of $P_{2,6}$. But then $P_{5,6}\cup P_{3,5}\cup P_{3,6}$ is an $F$-cycle, 
a contradiction. This proves (3).
\\
\\
(4) {\em There is no edge in $E(G)\setminus E(H)$ incident with $p_4$.}
\\
\\
Suppose that $p_4v\in E(G)\setminus E(H)$. From the symmetry, we may assume that $v$ belongs to the interior of $P_{1,5}\cup P_{5,6}\cup P_{2,6}$. But then adding the edges $p_4v$ and $p_1u$ to $H$, and deleting the interior of $P_{7,8}$ and the edge $p_3p_5$, gives a nonplanar graph, and so $E(P_{7,8})\subseteq F$. But $p_3p_7,p_3p_8\in F$, and so there is an $F$-cycle, a contradiction. This proves (4). 
\\
\\
(5) {\em If $u\ne p_2, p_3$, then $G$ is a conch.}
\\
\\
To show this, we must check that $p_3$
has no neighbour in the interiors of $P_{1,4}, P_{2,4}$, and there are no edges in $E(G)\setminus E(H)$ incident with $p_2$ except possibly $p_1p_2, p_2p_3$. Since $u\ne p_2, p_3$, it follows that $P_{1,4}$ has length one and the 
edge $p_1p_4\in F$. So $p_3$ has no neighbour in the interior of $P_{1,4}$. Suppose it has a neighbour $v$ in the interior of $P_{2,4}$.
Thus $P_{3,4}$ has length one and $p_3p_4, p_3v\in F$. Since $u \in V(P_{3,4}\cup P_{2,4})$, it follows that $u\in V(P_{2,4})$. 
Let $Q$ be the minimal subpath of $P_{2,4}$ from $p_4$ to $\{u,v\}$. Then $E(Q)\subseteq F$, 
and so there is an $F$-cycle (because $(p_1p_4,p_1u, p_3p_4,p_3v$ are all in $F$), a contradiction. So 
$p_3$
has no neighbour in the interiors of $P_{1,4}, P_{2,4}$. Next suppose that $p_2v\in E(G)\setminus E(H)$ where $v\ne p_1, p_3$.
By (3) (with $p_1,p_2$ exchanged) it follows that $v\in V(P_{3,4}\cup P_{1,4})$. Moreover, $P_{1,4}, P_{2,4}$ have length one; 
so both $u,v$ belong to $P_{3,4}$. Let $Q$ be the minimal path of $P_{3,4}$ between $p_4$ and $\{u,v\}$. Then $E(Q)\subseteq F$, and so there is
an $F$-cycle, since $p_1p_4,p_1u,p_2p_4,p_2v\in F$, a contradiction. This proves (5). 

\bigskip

From (5) we may assume that no edge in $E(G)\setminus E(H)$ is incident with $p_1$ except possible $p_1p_2,p_1p_3$, and similarly
none are incident with $p_2$ except possibly $p_1p_2,p_2p_3$. But then $G$ is a mussel. This proves \ref{conch}.~\bbox

In summary, then, by combining \ref{number2}, \ref{number4}, \ref{noU8}, \ref{doublewheel}, and \ref{conch}, we have proved the following, which implies \ref{planar} in view of \ref{reducetograph}:
\begin{thm}\label{summary}
Let $G$ be 3-connected and nonplanar. Then it is almost-planar if and only if $G$ is a M\"obius chain, or a double wheel,
or a conch, mussel, scallop, clam, or whelk.
\end{thm}


\section{Back to digraphs}\label{sec:backto}

\ref{summary} is satisfying, but does not really answer the original question:
what are the Kuratowski digraphs?
By \ref{reducetograph} and \ref{summary}, every Kuratowski digraph is an orientation
of one of the graphs of \ref{summary}. But not all the graphs of \ref{summary} can be oriented to make Kuratowski digraphs, 
and some can be oriented to do so in many ways. In this section we explore this issue.

If $G$ is a digraph, then for each $X\subseteq V(G)$, let $\delta^+(X)$ be the set of all edges of $G$ with tail in $X$ and head in $V(G)\setminus X$,
and let $\delta^-(X)=\delta^+(V(G)\setminus X)$. Let $\delta(X)$ be the set of all edges between $X$ and $V(G)\setminus X$.
If $G$ is a Kuratowshi digraph, then with $F$ defined as usual, for each 
$e=uv\in F$, since $G\setminus e$ is nonplanar and hence not strong (from the minimality property of $G$), there exists 
$A\subseteq V(G)$ with $\delta^+(A)=\{e\}$. We call such a set $A$ a {\em back-cut} for $e$. 
If in addition $\delta^-(A)\cap F=\emptyset$, 
we say $A$ is {\em clean}.

\begin{thm}\label{cleancut}
Let $G$ be a Kuratowski digraph, and let $F$ be the set of edges $e$ of $G$ such that $G\setminus e$ is nonplanar. Then 
there is a clean back-cut for each $e\in F$.
\end{thm}
\Proof 
Choose $P\subseteq F$ maximal such that 
for each $e\in F$, there is a back-cut $A$ for $e$ with $P\cap \delta^-(A)=\emptyset$. We suppose for a contradiction 
that $P\ne F$;
let $f=xy\in F\setminus P$. Let $B$ be a back-cut for $f$ with $P\cap \delta^-(B)=\emptyset$. From the maximality of $P$,
there exists $e=uv\in F$ such that every back-cut for $e$ contains an edge in $P\cup \{f\}$. Let $A$ be a back-cut for $e$
with $P\cap \delta^-(A)=\emptyset$, and consequently with $f\in \delta^-(A)$. Hence $x\in B\setminus A$
and $y\in A\setminus B$. 

Suppose that $v \notin A\cup B$.  Then $u\in B\setminus A$, since $e\notin \delta^+(B)$; and so
$A\cup B$ is also a back-cut for $e$. Moreover, $\delta^-(A\cup B)\subseteq \delta^-(A)\cup \delta^-(B)$, and so contains no edges in $P$; and $e\notin \delta^-(A\cup B)$, contradicting that every back-cut for $e$ contains an edge in $P\cup \{f\}$.
This proves that $v\in A\cup B$, and so $v\in B\setminus A$ since $v\notin A$. 

Suppose that $u\in A\setminus B$. 
Then $\delta^+(A\cap B)=\emptyset$, and so $A\cap B=\emptyset$, since $G$ is strong; and $\delta^+(A\cup B)=\emptyset$, and so
$A\cup B=V(G)$, since $G$ is strong. Thus $B=V(G)\setminus A$. But $e,f$ are the only edges between $A,B$, since 
any other edge would belong to $\delta^+(A)\cup \delta^+(B)$, which is impossible. Hence $G^-$ is not 3-connected, contrary to 
\ref{reducetograph}. 

This proves that $u\in A\cap B$.
But then $A\cap B$ is a back-cut for $e$, and $\delta^-(A\cap B)$ is disjoint from $P\cup \{f\}$ as before, a contradiction.
Hence $P=F$. This proves \ref{cleancut}.~\bbox

Let $G$ be an almost-planar graph, and let us say an orientation $H$ of $G$ is {\em good} if it is a Kuratowski 
digraph.  As usual, let $F$ be the set of edges $e$ such that $G\setminus e$ is nonplanar. 
Let us say a cycle $D$ of $G$ is {\em fundamental} if it has only one edge in $E(G)\setminus F$.
By \ref{cleancut}, in every good orientation, every fundamental cycle becomes a directed cycle. 

The result \ref{cleancut} is most helpful when $F$ is the edge-set of a spanning tree $T$ of $G$. 
In that case, for each $e\in F$, every clean back-cut for $e$ must be the vertex set of one of the two components of $T\setminus e$;
and so an orientation $H$ of $G$ is good if and only if for every $e\in E(G)\setminus F$, the fundamental cycle
containing $e$ is a directed cycle of $H$. 
It is easy to check that for any spanning tree in any 2-connected graph, there are at most two 
orientations of $G$ to make all the fundamental cycles directed, and they are reverses of each other. (Once we fix the direction of one edge, we know the direction of all edges in fundamental cycles containing this edge, and so on.)
Deciding whether there is a good orientation reduces in this case to a
2-SAT problem in which all clauses are equivalences, and this can be used to show that 
there is no good orientation if and only if for some odd $k\ge 3$, there are $k+1$ vertex-disjoint subtrees 
$T_0,T_1\LL T_k$ of $T$, such that for $1\le i\le k$ there is an edge of $T$ between $V(T_i)$ and $V(T_0)$, and for $1\le i\le k$
there is an edge of $E(G)\setminus E(T)$ between $V(T_i)$ and $V(T_{i+1})$ (reading subscripts modulo $k$). We omit the proof.

So for instance, the double wheel in Figure \ref{fig:nonplanar} has no good orientation. To see this, let $a,b$ be the two 
wheel centres, and let $C$ be the base cycle. Since  $a,b$ have a common neighbour in $C$, $F$ is indeed  the edge-set of a 
spanning tree. Moreover, there is a path $P$ of $C$ with both ends adjacent to $a$ and with a positive even number of its 
vertices adjacent to $b$,
and it is straightforward to check that this already precludes there being a good orientation. Let us say a double wheel is {\em proper}
if $G\setminus e$ is nonplanar for every edge $e$ between a wheel centre and the base cycle.  For proper double wheels in which the wheel centres have a common neighbour, one can show that it is necessary and sufficient
that there is no path $P$ as above. 

There are other types of almost-planar graphs in which $F$ makes a spanning tree; for instance, a mussel in which the hinge is adjacent to both tips.
Then one can show that there is a good orientation if and only if the three paths of the mussel from tip to tip that avoid 
the hinge either all have odd length or all have even length.

However, when $F$ does not make a spanning tree, it becomes harder to analyze when a good orientation exists, although it 
seems to become more likely that one does in fact exist.  Let $T$ be the forest with vertex set $V(G)$ and edge set $F$.
There are instances in which $T$ has two components and no good orientation exists; for instance, if $G$ is a proper double wheel,
and the two wheel centres have no common neighbour, then a good orientation exists if and only
if the base cycle has even length.
(Again, we omit the proof, which is another application of \ref{cleancut}.)

Let $e\in F$, and let $J$ be the minimal subgraph of $G$ containing $e$
such that every fundamental cycle with an edge in $J$ is a subgraph of $J$. It follows that the edges of $J$ in $F$
make a spanning tree of $J$; and $J$ is an induced subgraph of $G$. We call such a graph $J$ a {\em cell} of $G$. Every two
cells share at most one vertex, and so are edge-disjoint. (We skip these proofs, which are all easy.)
For a cell $J$, there are only two orientations of $J$ in which each fundamental cycle included in $J$ becomes directed, reverses of each other; and 
$H$ is a good orientation of $G$, and $J$ is a cell of $G$, then the restriction of $H$ to $J$ must be one of these two 
orientations. 

Thus, for a conch $G$ say, to decide whether it admits a good orientation, we need to know only a few things. 
Let $P,Q,R,p_1,p_2,p_3$ be as in the definition of a conch. We need to know whether $p_1$ is adjacent to $p_2$ or to $p_3$, and 
we need to know the parities of the lengths of $P,Q,R$. Nothing else is relevant; this information is enough  to
determine whether a good orientation exists, and we can easily decide each of the $2^5$ cases this yields. Similar, for each of the five types of graphs in Figure \ref{fig:clam}, the problem breaks into a small number of cases, each of which can be decided. 
(This is related to the ``fans'' mentioned in the introduction.) We have not bothered to figure out each case.

There remain M\"obius chains, and these are much more intractable (and much more interesting). 
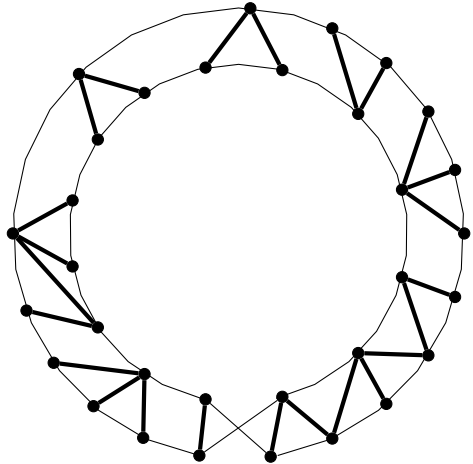
\begin{figure}[H]
\centering

\begin{tikzpicture}[scale=1.5,auto=left]
\tikzstyle{every node}=[inner sep=1.5pt, fill=black,circle,draw]
\def\r{1.5}
\def\a{360/12}

\node (a1) at ({\r*cos(9.5*\a)}, {\r*sin(9.5*\a)}){};
\node (a2) at ({\r*cos(10.5*\a)}, {\r*sin(10.5*\a)}){};
\node (a3) at ({\r*cos(11.5*\a)}, {\r*sin(11.5*\a)}){};
\node (a4) at ({\r*cos(.5*\a)}, {\r*sin(.5*\a)}){};
\node (a5) at ({\r*cos(1.5*\a)}, {\r*sin(1.5*\a)}){};
\node (a6) at ({\r*cos(2.5*\a)}, {\r*sin(2.5*\a)}){};

\def\a{360/16}
\node (a7) at ({\r*cos(4.5*\a)}, {\r*sin(4.5*\a)}){};
\node (a8) at ({\r*cos(5.5*\a)}, {\r*sin(5.5*\a)}){};
\node (a9) at ({\r*cos(6.5*\a)}, {\r*sin(6.5*\a)}){};
\node (a10) at ({\r*cos(7.5*\a)}, {\r*sin(7.5*\a)}){};
\node (a11) at ({\r*cos(8.5*\a)}, {\r*sin(8.5*\a)}){};
\node (a12) at ({\r*cos(9.5*\a)}, {\r*sin(9.5*\a)}){};
\node (a13) at ({\r*cos(10.5*\a)}, {\r*sin(10.5*\a)}){};
\node (a14) at ({\r*cos(11.5*\a)}, {\r*sin(11.5*\a)}){};
\draw [domain={-3*360/14:10*360/14}] plot ({\r*cos(\x)}, {\r*sin(\x)});

\def\r{2}
\def\a{360/22}
\node (b1) at ({\r*cos(17*\a)}, {\r*sin(17*\a)}){};
\node (b2) at ({\r*cos(18*\a)}, {\r*sin(18*\a)}){};
\node (b3) at ({\r*cos(19*\a)}, {\r*sin(19*\a)}){};
\node (b4) at ({\r*cos(20*\a)}, {\r*sin(20*\a)}){};
\node (b5) at ({\r*cos(21*\a)}, {\r*sin(21*\a)}){};
\node (b6) at ({\r*cos(0*\a)}, {\r*sin(0*\a)}){};
\node (b7) at ({\r*cos(1*\a)}, {\r*sin(1*\a)}){};
\node (b8) at ({\r*cos(2*\a)}, {\r*sin(2*\a)}){};
\node (b9) at ({\r*cos(3*\a)}, {\r*sin(3*\a)}){};
\node (b10) at ({\r*cos(4*\a)}, {\r*sin(4*\a)}){};
\node (b11) at ({\r*cos(87)}, {\r*sin(87)}){};

\def\a{360/20}

\node (b12) at ({\r*cos(135)}, {\r*sin(135)}){};
\node (b13) at ({\r*cos(180)}, {\r*sin(180)}){};
\node (b14) at ({\r*cos(200)}, {\r*sin(200)}){};
\node (b15) at ({\r*cos(215)}, {\r*sin(215)}){};
\node (b16) at ({\r*cos(230)}, {\r*sin(230)}){};
\node (b17) at ({\r*cos(245)}, {\r*sin(245)}){};
\node (b18) at ({\r*cos(260)}, {\r*sin(260)}){};

\draw [domain={-4*360/18:13*360/18}] plot ({\r*cos(\x)}, {\r*sin(\x)});

\foreach \from/\to in {a1/b18,a14/b1}
\draw [-] (\from) -- (\to);
\foreach \from/\to in {a1/b1,a14/b18, a1/b2,a2/b2,a2/b3,a2/b4,a3/b4,a3/b5,a4/b6,a4/b7,a4/b8,a5/b9,a5/b10,a6/b11,a7/b11,a8/b12,a9/b12,a10/b13,a11/b13,a12/b13, a12/b14,
a13/b15,a13/b16,a13/b17}
\draw [ultra thick] (\from) -- (\to);
\end{tikzpicture}
\caption{A M\"obius chain that does not admit a good orientation. The edges of $F$ are thick.} \label{fig:band}
\end{figure}

We have been able to give a pair of necessary and sufficient parity-type conditions for when a M\"obius chain admits a 
good orientation. But they are complicated, and the proof that they are correct is long and messy (about eight pages), and in the end we decided not to inflict it 
all
on the reader. For instance, here is the simpler of the necessary conditions: that if $u,v,w\in V(G)$ are distinct, 
and $uv,vw\in F$, and $u,v,w$ all have degree at least two in $T$, then $v$ has even degree in $T$. The other
necessary condition involves sequences (of arbitrary length) of components of $T$, and we do not give it here (as a clue, 
the graph in Figure
\ref{fig:band} violates it). 
It follows from that result that in a M\"obius chain: 
\begin{itemize}
\item every M\"obius chain in which $T$ contains no four-vertex path admits a good orientation; and
\item if $T$ is a spanning tree, then there is a good orientation if and only if every vertex of the spine of $T$ has
even degree in $T$ (where the {\em spine} is the unique minimal path of $T$ that meets all edges of $T$).
\end{itemize}

\end{document}